\documentclass[]{amsart}
\usepackage{url}
\usepackage{amssymb}
\usepackage{amsmath}
\usepackage{amsaddr}
\usepackage{stmaryrd}
\usepackage{siunitx}
\usepackage{commath}
\usepackage{enumitem}
\usepackage[caption=false]{subfig}
\usepackage{cleveref}
\usepackage{graphicx}
\usepackage{amsthm, amssymb, amsfonts, verbatim}
\makeatletter
\newcommand{\tnorm}{\@ifstar\@tnorms\@tnorm}
\newcommand{\@tnorm}[2][]{%
  \mathopen{#1|\mkern-1.5mu#1|\mkern-1.5mu#1|}
  #2
  \mathclose{#1|\mkern-1.5mu#1|\mkern-1.5mu#1|}
}
\makeatother

\newcommand{\jump}[1]{\llbracket #1 \rrbracket}
\newcommand{\av}[1]{\{\!\!\{#1\}\!\!\}}

\newtheorem{theorem}{Theorem}
\newtheorem{lemma}{Lemma}
\newtheorem{corollary}{Corollary}

\newcommand{\authorfootnotes}{\renewcommand\thefootnote{\@fnsymbol\c@footnote}}%

\begin{document}

\title[HDG for Stokes-Biot problem]{Hybridizable discontinuous Galerkin methods for the coupled Stokes--Biot problem}

\author{Aycil Cesmelioglu$^1$, Jeonghun J. Lee$^2$, Sander Rhebergen$^3$ }
\address{$^1$Department of Mathematics and Statistics, Oakland
  University, Michigan, USA \\ 
  $^2$ Department of Mathematics, Baylor University, Waco, Texas, USA \\
  $^3$ Department of Applied Mathematics, University of
  Waterloo, Ontario, Canada }

\email{$^1$cesmelio@oakland.edu, $^2$jeonghun\_lee@baylor.edu, $^3$srheberg@uwaterloo.ca}
 

\begin{abstract}
  We present and analyze a hybridizable discontinuous Galerkin (HDG)
  finite element method for the coupled Stokes--Biot problem. Of
  particular interest is that the discrete velocities and displacement
  are $H(\text{div})$-conforming and satisfy the compressibility
  equations pointwise on the elements. Furthermore, in the
  incompressible limit, the discretization is strongly
  conservative. We prove well-posedness of the discretization and,
  after combining the HDG method with backward Euler time stepping,
  present a priori error estimates that demonstrate that the method is
  free of volumetric locking. Numerical examples further demonstrate
  optimal rates of convergence in the $L^2$-norm for all unknowns and
  that the discretization is locking-free.
\end{abstract}


\date{July, 2022}
\maketitle

\section{Introduction}
\label{sec:introduction}

Many applications in environmental and biomedical engineering require
the modeling of the interaction between a free fluid and a deformable
porous media that is saturated by fluid. Such problems can be modeled
by the coupled Stokes--Biot problem, a model first proposed by
\cite{Showalter:2005}, in which the Stokes equations are coupled to
Biot's consolidation model \cite{Biot:1941,Biot:1955,Biot:1957}. At
the interface the two flows are coupled by enforcing mass
conservation, continuity of the normal stress, a balance of forces
across the interface, and the Beavers--Joseph--Saffman interface condition 
\cite{Beavers:1967,Saffman:1971}. Existence and uniqueness of weak
solutions to this problem were proven in \cite{Cesmelioglu:2017a}.

Various finite element methods have been proposed for the Stokes--Biot
problem. To better describe these methods we use $u^s$ and $p^s$ to
denote the fluid velocity and pressure in the Stokes region and $u^b$,
$z$, and $p^p$ to denote the displacement, Darcy velocity, and pore
pressure in the Biot region. \cite{Badia:2009} introduced the first
finite element method for this problem considering Lagrange finite
elements for $(u^s, p^s)$ and Lagrange and $H(\text{div})$ mixed
methods for $(u^b, z, p^p)$. Residual-based stabilization techniques
were used in both subdomains to stabilize this conforming finite
element method. A Lagrange multiplier method to weakly impose mass
conservation across the interface was studied in
\cite{Ambartsumyan:2018}. They used conforming stable methods for the
Stokes unknowns $(u^s, p^s)$
\cite{Brezzi:book,Arnold:1984,Taylor:1973,Girault:1986}, a Lagrange
finite element for $u^b$ and stable mixed finite elements for
$(z, p^p)$ \cite{Brezzi:book,Brezzi:1985,Raviart:1977}. To avoid
poroelastic locking when using a conforming stable finite element
method of the Stokes--Biot problem \cite{Cesmelioglu:2020a} studied
the use of a heuristic stabilization technique. A numerical study with
similar methods was done in \cite{Bergkamp:2020} for a Stokes--Biot
model with fluid entrance resistance. A finite element method based on
the total pressure formulation of the Biot model
\cite{JJLee:2017,Oyarzua:2016} was studied by
\cite{Baier:2022}. Eliminating the Darcy velocity $z$, they consider
stable Stokes elements for $(u^s, p^s)$ and $(u^b, p^b)$, where $p^b$
is the total pressure, and a piecewise continuous and polynomial space
for $p^p$. Finally, let us mention that a stress tensor based approach
using a weakly symmetric mixed method for poroelasticity \cite{Lee:2016}
and a conforming stable mixed method for Stokes was studied in
\cite{Ambartsumyan:2019,Li:2022}, and that partitioned time
discretization methods, focusing on efficient time stepping and
stability of partitioned schemes, are studied in
\cite{Bukac:2015,Bukac:2016,Oyekole:2020}.

In this paper we propose a hybridizable discontinuous Galerkin (HDG)
method for the coupled Stokes--Biot problem. The HDG method uses
hybridization to improve the computational efficiency of traditional
discontinuous Galerkin (DG) methods \cite{Cockburn:2009a}:
degrees-of-freedom (dofs) local to an element are eliminated resulting
in a global problem for dofs associated only with mesh facet
unknowns. The HDG method we propose here combines the HDG method of
\cite{Rhebergen:2017,Rhebergen:2018a} for the Stokes equations and the
HDG method studied in \cite{Cesmelioglu:2022} for the Biot model. (See
also \cite{Cesmelioglu:2020} for a similar HDG method and
\cite{Fu:2018} for an alternative approach for the coupled
Stokes--Darcy problem.) Our HDG method is constructed such that the discrete fluid velocities and displacement are H(div)-conforming. In addition, unlike existing finite element methods for the Stokes–Biot problem in the literature, the (semi-)discrete solution satisfies the compressibility equations and the mass equation in the poroelastic domain pointwise on each element of the mesh (where the latter holds provided the source/sink term lies in the pore pressure space).
These properties imply that the discretization is strongly
conservative \cite{Kanschat:2010} in the incompressible limit. In
addition to these exact approximation properties, we present an a
priori error analysis for the time-dependent coupled Stokes--Biot
problem that avoids using Gr\"onwall's inequality. A consequence of
the latter is that the space-time error estimates do not grow
exponentially in time.

The remainder of this paper is organized as follows. In
\cref{sec:stokesbiot,sec:hdg} we present the coupled Stokes--Biot
model and its HDG discretization. We proceed by proving consistency
and well-posedness of the HDG method in \cref{sec:consiswphdg} and
present an a priori error analysis of the discretization in
\cref{sec:apriorierrorhdg}. Numerical experiments in
\cref{sec:numerical_examples} verify our analysis and conclusions are
drawn in \cref{sec:conclusions}.

\section{The coupled Stokes--Biot equations}
\label{sec:stokesbiot}

Let $\Omega \subset \mathbb{R}^{d}$ be a polygonal (if $d=2$) or
polyhedral (if $d=3$) domain that is partitioned into two
non-overlapping domains $\Omega^s$ and $\Omega^b$. The Stokes
equations on $\Omega^s$ describe flow of an incompressible fluid while
a deformable porous structure on $\Omega^b$ is modelled by a
quasi-static poroelasticity model. We will assume that the interface
between the two subdomains,
$\Gamma^I = \overline{\partial\Omega}^s \cap
\overline{\partial\Omega^b}$ is polygonal. We denote by $n$ and $n^j$
the unit outward normal to, respectively, $\Omega$ and $\Omega^j$
($j=s,b$). On the interface $n = n^s = -n^b$. Furthermore, let us
denote by $J = (0, T]$ the time interval of interest.

In the poroelastic part of the domain we consider two different
partitions of the boundary
$\Gamma^b = \partial\Omega \cap \partial\Omega^b$. The first partition
is $\Gamma^b = \Gamma^b_P \cup \Gamma^b_F$ with
$\Gamma^b_P \cap \Gamma^b_F = \emptyset$ and $|\Gamma^b_P|>0$, while
the second partition is given by
$\Gamma^b = \Gamma^b_D \cup \Gamma^b_N$ where
$\Gamma^b_D \cap \Gamma^b_N = \emptyset$, $|\Gamma^b_D|>0$, and $|\Gamma^b_N|>0$. In the
fluid domain we partition the boundary
$\Gamma^s = \partial\Omega \cap \partial\Omega^s$ as
$\Gamma^s = \Gamma^s_D \cup \Gamma^s_N$ with
$\Gamma^s_D \cap \Gamma^s_N = \emptyset$, $|\Gamma^s_D|>0$, and $|\Gamma_N^s|>0$.

We denote body forces by $f^s :\Omega^s \to \mathbb{R}^{d}$ and
$f^b :\Omega^b \to \mathbb{R}^{d}$ and the source/sink term by
$g^b :\Omega^b \to \mathbb{R}$. Furthermore, we denote by $\mu^s > 0$
the (constant) dynamic viscosity of the fluid, the Biot--Willis
constant and the specific storage coefficient are denoted by,
respectively, $0 < \alpha < 1$ and $c_0 \ge 0$, the positive
permeability constant is denoted by $\kappa$, and the Lam\'e constants
are denoted by $\mu^b$ and $\lambda$. Note that Young's modulus of
elasticity $E$ and Poisson's ratio $\nu$ are related to the Lam\'e
constant by $E = (3\lambda+2\mu^b)\mu^b/(\lambda+\mu^b)$ and
$\nu = \lambda/(2(\lambda + \mu^b))$.

Using the total pressure formulation of Biot's consolidation model
\cite{JJLee:2017,Oyarzua:2016}, we can then state the coupled
Stokes--Biot problem as: Find the fluid velocity
$u^s:\Omega^s \to \mathbb{R}^{d}$, the fluid pressure
$p^s :\Omega^s \to \mathbb{R}$, the solid displacement
$u^b:\Omega^b \to \mathbb{R}^{d}$, the pore pressure
$p^p:\Omega^b \to \mathbb{R}$, the total pressure
$p^b:\Omega^b \to \mathbb{R}$, and the Darcy velocity
$z:\Omega^b \to \mathbb{R}^{d}$ such that
\begin{subequations}
  \label{eq:stokesbiot}
  \begin{align}
    \label{eq:stokesbiot_1}
    -\nabla \cdot \sigma^j &= f^j && \text{in } \Omega^j\times J, \quad j=s,b,
    \\
    \label{eq:stokesbiot_2a}
    -\nabla \cdot u^s 
    &= 0 && \text{in } \Omega^s\times J,
    \\
    \label{eq:stokesbiot_2b}
    -\nabla \cdot u^b + \lambda^{-1}(\alpha p^p - p^b)
    &= 0 && \text{in } \Omega^b\times J,
    \\
    \label{eq:stokesbiot_3}
    c_0 \partial_t p^p + \alpha \lambda^{-1} (\alpha \partial_t p^p - \partial_t p^b) + \nabla \cdot z
                           &= g^b && \text{in }\Omega^b \times J,
    \\
    \label{eq:stokesbiot_4}
    \kappa^{-1} z + \nabla p^p
    &= 0 && \text{in }\Omega^b \times J,
  \end{align}  
\end{subequations}
and such that
\begin{subequations}
  \label{eq:interface}
  \begin{align}
    \label{eq:bc_I_u}
    u^s\cdot n &= \del[1]{\partial_t u^b + z}\cdot n & & \text{on } \Gamma_I \times J,
    \\
    \label{eq:bc_I_ss_sb}
    \sigma^sn &= \sigma^bn & & \text{on } \Gamma_I \times J,
    \\
    \label{eq:bc_I_p}
    -(\sigma^s n)\cdot n &= p^p & & \text{on } \Gamma_I \times J,
    \\    
    \label{eq:bc_I_slip}
    -2\mu^s\del[0]{\varepsilon(u^s)n}^t &= \gamma (\mu^s/\kappa)^{1/2}(u^s-\partial_t u^b)^t
                                                     & & \text{on } \Gamma_I \times J.
  \end{align}
\end{subequations}
Here $\sigma^j := 2\mu^j\varepsilon(u^j) - p^j\mathbb{I}$ ($j=s,b$),
$\varepsilon(u) := \del[0]{\nabla u + (\nabla u)^T}/2$,
$(w)^t:= w - (w\cdot n)n$, and \cref{eq:bc_I_slip} is the
Beavers--Joseph--Saffman condition \cite{Beavers:1967,Saffman:1971} in
which $\gamma > 0$ is an experimentally determined dimensionless
constant. To close the system we impose the boundary conditions
\begin{subequations}
  \label{eq:bics}
  \begin{align}
    \label{eq:bics-us}
    u^j &= 0 && \text{on } \Gamma^j_D\times J, \quad j=s,b,
    \\
    \label{eq:bics-sigmas}
    \sigma^j n &= 0 &&  \text{on } \Gamma^j_N\times J, \quad j=s,b,
    \\    
    \label{eq:bics-pb}
    p^p &= 0 && \text{on } \Gamma^b_P\times J,
    \\
    \label{eq:bics-zb}
    z \cdot n &= 0 && \text{on } \Gamma^b_F\times J,
  \end{align}
\end{subequations}
and initial conditions $p^p(x,0) =p_0^p(x)$ in $\Omega^b$ and
$u^b(x,0) = u^b_0(x)$ in $\Omega^b$.

For notational convenience it will be useful to define functions $u$
and $p$ on the whole domain $\Omega$ which are such that
$u|_{\Omega^j} = u^j$ and $p|_{\Omega^j}=p^j$ for $j=s,b$.

In the a priori error analysis in \cref{sec:apriorierrorhdg} we will
assume that there exists a $\nu_*$ such that $0 < \nu_* \le \nu < 0.5$
on $\Omega^b$ implying that $C_* \mu^b \le \lambda$ with
$C_* = 2\nu_*/(1-2\nu_*)$.

\section{Discretizing the Stokes--Biot problem}
\label{sec:hdg}

\subsection{Notation}
\label{ss:notation}

Let $\mathcal{T}^j$ denote a shape-regular triangulation of
$\Omega^j$, $j=s, b$. We will assume that $\mathcal{T}^s$ and
$\mathcal{T}^b$ match at the interface and define
$\mathcal{T}:=\mathcal{T}^s\cup \mathcal{T}^b$. For any element
$K\in \mathcal{T}$, $h_K$ denotes its diameter and
$h:=\max_{K\in \mathcal{T}} h_K$ defines the meshsize of the
triangulation. We denote the set of all facets by $\mathcal{F}$, the
set of all facets in $\bar{\Omega}^j$ by $\mathcal{F}^j$, $j=s, b$,
the set of all facets in $\Omega^j$ by $\mathcal{F}_{int}^j$, and the
set of all facets that lie on $\Gamma_I$, $\Gamma_N^j$, $\Gamma_F^b$,
and $\Gamma_P^b$ by $\mathcal{F}_I$, $\mathcal{F}_N^j$,
$\mathcal{F}_F^b$, and $\mathcal{F}_P^b$, respectively. Finally, we
set $\Gamma_0^j=\cup_{F\in \mathcal{F}^j}F$, $j=s,b$.

Various approximation spaces are required to define the HDG
discretization of the Stokes--Biot problem
\cref{eq:stokesbiot,eq:interface,eq:bics}. These approximation spaces
are discontinuous Galerkin (DG) spaces defined on $\Omega$ and
$\Omega^j$ ($j=s,b$):
\begin{subequations}
  \label{eqn:spaces_cell}
  \begin{align}
    \label{eqn:spaces_cell_Vj}
    V_h^j
    &:= \cbr[1]{v_h\in \sbr[0]{L^2(\Omega^j)}^d
      : \ v_h \in \sbr[0]{P_k(K)}^d, \ \forall\ K\in\mathcal{T}^j}, \quad j=s,b,
    \\
    \label{eqn:spaces_cell_Qj}
    Q_h^j
    &:= \cbr[1]{q_h\in L^2(\Omega^j) : \ q_h \in P_{k-1}(K) ,\
      \forall \ K \in \mathcal{T}^j},  \quad j=s,b,
    \\
    \label{eqn:spaces_cell_V}
    V_h
    &:= \cbr[1]{v_h \in \sbr[0]{L^2(\Omega)}^d : \ v_h \in \sbr[0]{P_k(K)}^d,\ \forall K \in \mathcal{T}},
    \\
    \label{eqn:spaces_cell_Q}
    Q_h 
    &:= \cbr[1]{q_h\in L^2(\Omega) : \ q_h \in P_{k-1}(K) ,\
      \forall \ K \in \mathcal{T}}.
  \end{align}
\end{subequations}
Note that for functions $u_h \in V_h$ and $p_h \in Q_h$ we have,
respectively, that $u_h|_{\Omega^j} = u_h^j \in V_h^j$ and
$p_h|_{\Omega^j} = p_h^j \in Q_h^j$ for $j=s,b$. The HDG
discretization also requires the following facet DG spaces that are
defined on $\Gamma^j_0$ ($j=s,b$):
\begin{subequations}
  \begin{align}
    \label{eqn:HDG_spaces_facet_Vbar}
    \bar{V}_h^j
    &:= \cbr[1]{\bar{v}_h \in \sbr[0]{L^2(\Gamma_0^j)}^d:\ \bar{v}_h \in
      \sbr[0]{P_{k}(F)}^d\ \forall\ F \in \mathcal{F}^j,\ \bar{v}_h
      = 0 \text{ on } \Gamma_D^j},
    \\
    \label{eqn:HDG_spaces_facet_Qbar}
    \bar{Q}_h^j
    &:= \cbr[1]{\bar{q}_h \in L^2(\Gamma_0^j) : \ \bar{q}_h \in
      P_{k}(F) \ \forall\ F \in \mathcal{F}^j}, \quad j=s,b,
    \\
    \label{eqn:HDG_spaces_facet_Qbar0}
    \bar{Q}_h^{b0}
    &:=\cbr[0]{\bar{q}_h \in \bar{Q}_h^b:
      \bar{q}_h=0 \text{ on } \Gamma_P^b}.
  \end{align}
  \label{eqn:HDG_spaces_facet}
\end{subequations}
For notational purposes, we group cell and facet unknowns as follows:
\begin{align*}
  \boldsymbol{v}_h = (v_h, \bar{v}_h^s, \bar{v}_h^b) &\in \boldsymbol{V}_h := V_h \times \bar{V}_h^s \times \bar{V}_h^b, &
  \\
  \boldsymbol{q}_h = (q_h, \bar{q}_h^s, \bar{q}_h^b) &\in \boldsymbol{Q}_h := Q_h \times \bar{Q}_h^s \times \bar{Q}_h^b,
  \\
  \boldsymbol{q}_h^p = (q_h^p, \bar{q}_h^p) &\in \boldsymbol{Q}_h^{b0} := Q_h^b \times \bar{Q}_h^{b0},
  \\
  (\boldsymbol{v}_h, \boldsymbol{q}_h, w_h, \boldsymbol{q}_h^p) &\in \boldsymbol{X}_h
  := \boldsymbol{V}_h \times \boldsymbol{Q}_h \times V_h^b \times
  \boldsymbol{Q}_h^{b0},
\end{align*}
and for $j=s,b$, 
\begin{equation*}
  \boldsymbol{v}_h^j = (v_h^j, \bar{v}_h^j) \in \boldsymbol{V}_h^j := V_h^j \times \bar{V}_h^j,
  \qquad
  \boldsymbol{q}_h^j = (q_h^j, \bar{q}_h^j) \in \boldsymbol{Q}_h^j := Q_h^j \times \bar{Q}_h^j.  
\end{equation*}

We will use standard innerproduct notation: for scalar functions $p$
and $q$ on an element $K$, $(p,q)_K = \int_{K}pq\dif x$; for vector
functions $p$ and $q$, $(p,q)_K = \int_{K} p \cdot q \dif x$; and for
matrix functions $p$ and $q$, $(p,q)_K = \int_K p : q \dif x$. Let $D$
be the $(d-1)$-dimensional boundary $\partial K$ or facet
$F \subset \partial K$ of an element $K$. We then write
$\langle p, q \rangle_D = \int_D p \odot q \dif s$, where $\odot$ is
multiplication if $p$ and $q$ are scalar functions, the dot-product if
$p$ and $q$ are vector functions, and the dyadic product if $p$ and
$q$ are matrix functions. We furthermore define
$(p,q)_{\Omega^j} := \sum_{K \in \mathcal{T}^j} (p,q)_K$ and
$\langle p,q \rangle_{\partial\mathcal{T}^j} := \sum_{K \in
  \mathcal{T}^j}\langle p, q \rangle_{\partial K}$ for $j=s,b$ and
$(p,q)_{\Omega} := \sum_{K \in \mathcal{T}} (p, q)_K$ and
$\langle p, q \rangle_{\partial \mathcal{T}} := \sum_{K \in
  \mathcal{T}} \langle p, q \rangle_{\partial K}$, while on the
interface $\Gamma_I$ we define
$\langle p, q \rangle_{\Gamma_I} := \sum_{F \in \mathcal{F}_I} \langle
p, q \rangle_F$. At this point it will be useful to also define
$\norm{q}_{K} := (q,q)_{K}^{1/2}$,
$\norm{q}_{\partial K} := \langle q, q\rangle_{\partial K}^{1/2}$,
$\norm{q}_F := \langle q, q \rangle_{F}^{1/2}$,
$\norm{q}_{\Gamma_I} := \langle q, q \rangle_{\Gamma_I}^{1/2}$, and
$\norm{q}_{\Omega^j} := (q,q)_{\Omega^j}^{1/2}$ for $j=s,b$.

The following bilinear forms are used in the following sections to
define the HDG method (where $j=s,b$):
\begin{subequations}
  \begin{align*}
    a_h^j(\boldsymbol{u}^j, \boldsymbol{v}^j)
    :=&
        (2\mu^j\varepsilon(u), \varepsilon(v))_{\Omega^j}
        + \sum_{K\in\mathcal{T}^j}\langle 2\beta^j\mu^j h_K^{-1} (u-\bar{u}^j), (v-\bar{v}^j) \rangle_{\partial K}
    \\
    \nonumber    
      & - \langle 2\mu^j \varepsilon(u)n^j, (v-\bar{v}^j) \rangle_{\partial\mathcal{T}^j}
        - \langle 2\mu^j \varepsilon(v)n^j, (u-\bar{u}^j) \rangle_{\partial\mathcal{T}^j},
    \\
    \nonumber
    a_h(\boldsymbol{u}, \boldsymbol{v})
    :=& a_h^s(\boldsymbol{u}^s, \boldsymbol{v}^s) + a_h^b(\boldsymbol{u}^b, \boldsymbol{v}^b),
    \\
    \nonumber
    b_h^{j}(\boldsymbol{q}^j, \boldsymbol{v}^j)
    :=&
        -(q, \nabla \cdot v)_{\Omega^j} + \langle \bar{q}^j, (v-\bar{v}^j) \cdot n^j \rangle_{\partial\mathcal{T}^j},
    \\
    \nonumber
    b_h(\boldsymbol{q}, \boldsymbol{v})
    :=& b_h^{s}(\boldsymbol{q}^s, \boldsymbol{v}^s) + b_h^{b}(\boldsymbol{q}^b, \boldsymbol{v}^b).
    \\
    \nonumber
    c_h((p,r),q)
    :=& (\lambda^{-1}(\alpha p - r), q)_{\Omega^b},
    \\
    \nonumber
    a_h^I((\bar{u}^s,\bar{u}^b),(\bar{v}^s,\bar{v}^b))
    :=& \langle \gamma(\mu^s/\kappa)^{1/2}(\bar{u}^s-\bar{u}^b)^t, (\bar{v}^s-\bar{v}^b)^t \rangle_{\Gamma_I},
    \\
    \nonumber
    b_h^I(\bar{q}^p,(\bar{v}^s,\bar{v}^b))
    :=& \langle \bar{q}^p, (\bar{v}^s-\bar{v}^b) \cdot n \rangle_{\Gamma_I}.
  \end{align*}
\end{subequations}

\subsection{The HDG method}
\label{ss:semihdg}

Let $f|_{\Omega^j} = f^j$ for $j=s,b$. We propose the following
semi-discrete HDG method for the coupled Stokes and Biot problem
\cref{eq:stokesbiot,eq:interface,eq:bics}: Find
$(\boldsymbol{u}_h(t), \boldsymbol{p}_h(t), z_h(t),
\boldsymbol{p}_h^p(t)) \in \boldsymbol{X}_h$, for $t\in J$,
such that
{\small
\begin{subequations}
  \label{eq:semidiscreteHDG}
  \begin{align}
    a_h(\boldsymbol{u}_h, \boldsymbol{v}_h) 
    + a_h^I((\bar{u}_h^s,\partial_t\bar{u}_h^b),(\bar{v}_h^s,\bar{v}_h^b))
    + b_h(\boldsymbol{p}_h, \boldsymbol{v}_h)
    + b_h^I(\bar{p}_h^p,(\bar{v}_h^s,\bar{v}_h^b))
    =& (f, v_h)_{\Omega},
    \\
    b_h(\boldsymbol{q}_h, \boldsymbol{u}_h)
    + c_h((p_h^p,p_h^b),q_h^b)
    =& 0, \label{eq:semidiscreteHDG-b}
    \\
    (c_0\partial_t p_h^p, q_h^p)_{\Omega^b}
    + c_h((\partial_tp_h^p,\partial_tp_h^b), \alpha q_h^p)
    - b_h^{b}(\boldsymbol{q}_h^p, (z_h, 0))
    - b_h^I(\bar{q}_h^p,(\bar{u}_h^s, \partial_t\bar{u}_h^b))
    =& (g^b, q_h^p)_{\Omega^b},
    \\
    (\kappa^{-1} z_h, w_h)_{\Omega^b} + b_h^{b}(\boldsymbol{p}_h^p, (w_h, 0))
    =& 0,  
  \end{align}  
\end{subequations}
}
for all
$(\boldsymbol{v}_h, \boldsymbol{q}_h, w_h, \boldsymbol{q}_h^p)\in
\boldsymbol{X}_h$.

Let us partition the time interval $J$ using the time levels
$0=t^0 < t^1 < \hdots < t^N = T$ where, for $n=0,\hdots,N$,
$t^n = n\Delta t$ with $\Delta t > 0$ the time step. Functions $g$ and
$g^j$ (with $j=s,b$) evaluated at time level $t^i$ are denoted by,
respectively, $g^i$ and $g^{j,i}$. We furthermore denote the first
order discrete time derivative as $d_tg^n := (g^n - g^{n-1})/\Delta t$
for $n=1,\hdots,N$. (It will be clear from context whether $n$ denotes
a time-level or a normal vector.) Applying Backward Euler
time-stepping to \cref{eq:semidiscreteHDG} we obtain the
fully-discrete method: Find
$(\boldsymbol{u}_h^{n+1}, \boldsymbol{p}_h^{n+1}, z_h^{n+1},
\boldsymbol{p}_h^{p,n+1})\in \boldsymbol{X}_h$ such that
\begin{subequations}
  \label{eq:fullydiscreteHDG}
  \begin{align}
    \label{eq:fullydiscreteHDG-a}
    &a_h(\boldsymbol{u}_h^{n+1}, \boldsymbol{v}_h) 
      + a_h^I((\bar{u}_h^{s,n+1},d_t\bar{u}_h^{b,n+1}),(\bar{v}_h^s,\bar{v}_h^b))
    \\ \nonumber
    &\hspace{7em}+ b_h(\boldsymbol{p}_h^{n+1}, \boldsymbol{v}_h)
      + b_h^I(\bar{p}_h^{p,n+1},(\bar{v}_h^s,\bar{v}_h^b))
      = (f^{n+1}, v_h)_{\Omega},
    \\
    \label{eq:fullydiscreteHDG-b}
    &b_h(\boldsymbol{q}_h, \boldsymbol{u}_h^{n+1})
      + c_h((p_h^{p,n+1},p_h^{b,n+1}),q_h^b)
      = 0,
    \\
    \label{eq:fullydiscreteHDG-c}
    &(c_0 d_t p_h^{p,n+1}, q_h^p)_{\Omega^b}
      + c_h((d_tp_h^{p,n+1},d_tp_h^{b,n+1}), \alpha q_h^p)
      - b_h^{b}(\boldsymbol{q}_h^p, (z_h^{n+1}, 0))
    \\ \nonumber 
    &\hspace{7em}- b_h^I(\bar{q}_h^p,(\bar{u}_h^{s,n+1}, d_t\bar{u}_h^{b,n+1}))
      = (g^{b,n+1}, q_h^p)_{\Omega^b},
    \\
    \label{eq:fullydiscreteHDG-d}
    &(\kappa^{-1} z_h^{n+1}, w_h)_{\Omega^b} + b_h^{b}(\boldsymbol{p}_h^{p,n+1}, (w_h, 0))
      = 0,  
  \end{align}  
\end{subequations}
for all
$(\boldsymbol{v}_h, \boldsymbol{q}_h, w_h, \boldsymbol{q}_h^p)\in
\boldsymbol{X}_h$.

\begin{lemma}
  \label{lem:properties-HDG-solution}
  The following properties of the solution to
  \cref{eq:fullydiscreteHDG} hold:
  \begin{subequations}
    \label{eq:uhdiv}
    \begin{align}
      \label{eq:uhdiv-a}
      \jump{u_h^{j,n+1} \cdot n} =& 0 && \forall x \in F, && \forall F \in \mathcal{F}_{int}^j \cup \mathcal{F}_D^j,
      \\
      \label{eq:uhdiv-b}
      u_h^{j,n+1} \cdot n =& \bar{u}_h^{j,n+1} \cdot n && \forall x \in F, && \forall F \in \mathcal{F}_N^j \cup \mathcal{F}_I,
      \\
      \label{eq:zhdiv-a}
      \jump{z_h^{n+1} \cdot n} =& 0 && \forall x \in F, && \forall F \in \mathcal{F}_{int}^b \cup \mathcal{F}_F^b,
      \\
      \label{eq:zhdiv-d}
      u_h^{s,n+1} \cdot n =& (z_h^{n+1} + d_tu_h^{b,n+1})\cdot n && \forall x \in F, && \forall F \in \mathcal{F}_I,
      \\
      \label{eq:uhdiv-c}
      \nabla \cdot u_h^{s,n+1} =& 0 && \forall x \in K, && \forall K \in \mathcal{T}^s,
      \\
      \label{eq:uhbdiv-e}
      \nabla \cdot u_h^{b,n+1} =& \lambda^{-1}(\alpha p_h^{p,n+1} - p_h^{b,n+1}) && \forall x \in K, && \forall K \in \mathcal{T}^b,
      \\
      \label{eq:zhdiv-f}
      \nabla \cdot z_h^{n+1} =& \Pi_Q^b g^b -c_0d_tp_h^{p,n+1} && &&
      \\ \nonumber
      & - \alpha\lambda^{-1}(\alpha d_tp_h^{p,n+1}-d_tp_h^{b,n+1}) && \forall x \in K, && \forall K \in \mathcal{T}^b,
    \end{align}
  \end{subequations}
  where $\Pi_Q^b$ is the $L^2$-projection operator onto $Q_h^b$ and
  where $\jump{\cdot}$ is the usual jump operator.
\end{lemma}
\begin{proof}
  Setting $q_h=0$ in \cref{eq:fullydiscreteHDG-b}, we find for all
  $\bar{q}_h^j \in \bar{Q}_h^j$ with $j=s,b$ that
  \begin{equation*}
    \begin{split}
      0
      =& \sum_{K \in \mathcal{T}^j}\langle \bar{q}_h^j, (u_h^{j,n+1} - \bar{u}_h^{j,n+1})\cdot n^j \rangle_{\partial K}
      \\
      =& \sum_{F \in \mathcal{F}^j_{int}\cup \mathcal{F}_D^j} \langle \bar{q}_h^j, \jump{u_h^{j,n+1}\cdot n^j} \rangle_F
      + \sum_{F \in (\mathcal{F}_N^j \cup \mathcal{F}_I)} \langle \bar{q}_h^j, (u_h^{j,n+1} - \bar{u}_h^{j,n+1})\cdot n^j \rangle_{F}. 
    \end{split}
  \end{equation*}
  \Cref{eq:uhdiv-a,eq:uhdiv-b} follow noting that
  $(u_h^{j,n+1}\cdot n)|_F \in P_k(F)$ and
  $(\bar{u}_h^{j,n+1}\cdot n)|_F \in P_k(F)$.
  Setting now $q_h^p = 0$ in \cref{eq:fullydiscreteHDG-c}, we find for
  all $\bar{q}_h^p \in \bar{Q}_h^{b0}$ that
  \begin{equation*}
    \begin{split}
      0
      =& \sum_{K \in \mathcal{T}^b}\langle \bar{q}_h^p, z_h^{n+1} \cdot n^b \rangle_{\partial K}
      + \langle \bar{q}_h^p, (\bar{u}_h^{s,n+1} - d_t\bar{u}_h^{b,n+1})\cdot n \rangle_{\Gamma_I}
      \\
      =& \sum_{F \in \mathcal{F}^b_{int}} \langle \bar{q}_h^p, \jump{z_h^{n+1}\cdot n^b} \rangle_F    
      + \sum_{F \in \mathcal{F}_F^b} \langle \bar{q}_h^p, z_h^{n+1} \cdot n^b \rangle_{F}
       \\
       & + \langle \bar{q}_h^p, (\bar{u}_h^{s,n+1} - (z_h^{n+1} + d_t\bar{u}_h^{b,n+1}))\cdot n \rangle_{\Gamma_I}
    \end{split}
  \end{equation*}
  where the second equality is because $n=-n^b$.
  
  \Cref{eq:zhdiv-a} follows since $(z_h^{n+1}\cdot n)|_F \in P_k(F)$,
  \cref{eq:zhdiv-d} is an immediate consequence of \cref{eq:uhdiv-b}
  and since
  $((\bar{u}_h^{s,n+1} - (z_h^{n+1} + d_t\bar{u}_h^{b,n+1}))\cdot
  n)|_F \in P_k(F)$, while \cref{eq:uhdiv-c} follows after setting
  $\boldsymbol{q}_h^b = \boldsymbol{0}$ and $\bar{q}_h^s = 0$ in
  \cref{eq:fullydiscreteHDG-b} and noting that
  $\nabla \cdot u_h^s \in P_{k-1}(K)$. Setting
  $\boldsymbol{q}_h^s = \boldsymbol{0}$, $\bar{q}_h^b = 0$, and
  $q_h^b = -\nabla \cdot u_h^{b,n+1} + \lambda^{-1}(\alpha p_h^{p,n+1}
  - p_h^{b,n+1})$ in \cref{eq:fullydiscreteHDG-b} results in
  \cref{eq:uhbdiv-e}. Finally, \cref{eq:zhdiv-f} follows by setting
  $\bar{q}_h^p = 0$ and
  $q_h^p = c_0d_tp_h^{p,n+1} + \alpha\lambda^{-1}(\alpha
  d_tp_h^{p,n+1}-d_tp_h^{b,n+1}) + \nabla \cdot z_h^{n+1} - \Pi_Q^b
  g^b$ in \cref{eq:fullydiscreteHDG-c}. 
\end{proof}

\Cref{lem:properties-HDG-solution} demonstrates that $u_h^{s,n+1}$,
$u_h^{b,n+1}$, and $z_h^{n+1}$ are $H(\text{div})$-conforming and that
the compressibility equations \cref{eq:stokesbiot_2a,eq:stokesbiot_2b}
are satisfied pointwise on the elements by the numerical solution. For
the semi-discrete method \cref{eq:semidiscreteHDG}, \cref{eq:zhdiv-f}
can be replaced by
\begin{equation*}
  c_0 \partial_t p_h^p + \alpha\lambda^{-1}(\alpha\partial_tp_h^p - \partial_tp_h^b) + \nabla \cdot z_h
  = \Pi_Qg^b \qquad \forall x \in K, \forall K \in \mathcal{T}_h^b, \forall t \in J,
\end{equation*}
which states that mass is conserved pointwise on the elements if
$g^b \in Q_h^b$. In the incompressible limit, i.e.,
$\lambda \to \infty$ and $c_0 \to 0$, the HDG method is strongly
conservative \cite{Kanschat:2010}.

\section{Consistency, stability, existence and uniqueness}
\label{sec:consiswphdg}

This section is devoted to proving consistency and stability of the semi-discrete
HDG method \cref{eq:semidiscreteHDG} and existence and uniqueness of a
solution to the fully-discrete HDG method
\cref{eq:fullydiscreteHDG}. We start with consistency.

\begin{lemma}[Consistency]
  \label{lem:consistency}
  Let $(u,p,z,p^p)$ be a solution to the coupled Stokes--Biot problem
  \cref{eq:stokesbiot,eq:interface,eq:bics}. Let
  $\boldsymbol{u} := (u, u|_{\Gamma_0^s}, u|_{\Gamma_0^b})$,
  $\boldsymbol{p} := (p, p|_{\Gamma_0^s}, u|_{\Gamma_0^b})$, and
  $\boldsymbol{p}^p := (p^p, p^p|_{\Gamma_0^b})$. Then
  $(\boldsymbol{u}, \boldsymbol{p}, z, \boldsymbol{p}^p)$ satisifes
  the semi-discrete problem \cref{eq:semidiscreteHDG}.
\end{lemma}
\begin{proof}
  Integrating by parts and using the smoothness of $u^j$ and
  single-valuedness of $\bar{v}_h^j$ ($j=s,b$),
  \begin{equation*}
    \label{eq:consis_ahj0-1}
    \begin{split}
      a_h(\boldsymbol{u}, \boldsymbol{v}_h)
      =& \sum_{j=s,b}\big( -(2\mu^j(\nabla \cdot \varepsilon(u^j)), v_h^j )_{\Omega^j}
      +\langle 2\mu^j\varepsilon(u^j)n^j, \bar{v}_h^j\rangle_{\Gamma_I} \big) 
      \\
      &+ \langle 2\mu^s\varepsilon(u^s)n^s, \bar{v}_h^s\rangle_{\Gamma_N^s}
      + \langle 2\mu^b\varepsilon(u^b)n^b, \bar{v}_h^b\rangle_{\Gamma_N^b}.
    \end{split}
  \end{equation*}
  Similarly,
  \begin{equation*}
    \label{eq:bhscons-1}
    \begin{split}
      b_h(\boldsymbol{p}, \boldsymbol{v}_h)
      =&
      (\nabla p, v_h)_{\Omega}
      - \langle p^s, \bar{v}_h^s \cdot n^s\rangle_{\Gamma_I} - \langle p^b, \bar{v}_h^b \cdot n^b \rangle_{\Gamma_I}
       \\
       &
      - \langle p^s, \bar{v}_h^s \cdot n^s\rangle_{\Gamma_N^s} -  \langle p^b, \bar{v}_h^b \cdot n^b \rangle_{\Gamma_N^b}.
    \end{split}
  \end{equation*}
  Hence by \cref{eq:stokesbiot_1,eq:bics-sigmas}, 
   \begin{equation}
    \label{eq:consis_ahbh}
    \begin{split}
      a_h(\boldsymbol{u}, \boldsymbol{v}_h) + b_h(\boldsymbol{p}, \boldsymbol{v}_h)
      =&\langle 2\mu^s\varepsilon(u^s)n^s, \bar{v}_h^s\rangle_{\Gamma_I}
      + \langle 2\mu^b\varepsilon(u^b)n^b, \bar{v}_h^b\rangle_{\Gamma_I} \\
      &- \langle p^s, \bar{v}_h^s\cdot n^s \rangle_{\Gamma_I}
      - \langle p^b, \bar{v}_h^b\cdot n^b \rangle_{\Gamma_I} + (f^s,v_h^s)_{\Omega^s} + (f^b,v_h^b)_{\Omega^b}.
    \end{split}
   \end{equation}
   By \cref{eq:bc_I_u,eq:bc_I_ss_sb,eq:bc_I_p,eq:bc_I_slip}, on $\Gamma_I$
   \begin{equation}
     \label{eq:IC-j}
     \begin{split}
       2\mu^j\varepsilon(u^j)n^j
       &= (n^j\cdot 2\mu^j\varepsilon(u^j)n^j) n^j + (2\mu^j\varepsilon(u^j)n^j)^t\\
       &= p^j n^j - p^p n^j - \zeta^j \gamma (\mu^s/\kappa)^{1/2}(u^s-\partial_t u^b)^t,
     \end{split}
   \end{equation}
   where $\zeta^s = 1$ and $\zeta^b=-1$. Combining
   \cref{eq:consis_ahbh} with \cref{eq:IC-j} and using that
   $n = n^s = -n^b$,
   \begin{equation*}
    \begin{split}
      a_h(\boldsymbol{u}, \boldsymbol{v}_h) + b_h(\boldsymbol{p}, \boldsymbol{v}_h)
      =&-\langle\gamma (\mu^s/\kappa)^{1/2}(u^s-\partial_t u^b)^t,  (\bar{v}_h^s -\bar{v}_h^b)^t\rangle_{\Gamma_I}
      \\
      &- \langle p^p, (\bar{v}_h^s - \bar{v}_h^b) \cdot n^s\rangle_{\Gamma_I}  + (f,v_h)_{\Omega},
    \end{split}
   \end{equation*}
   which, after rearranging, results in
   \begin{equation}
     \label{eq:consis_ahbh-2}
     a_h(\boldsymbol{u}, \boldsymbol{v}_h) + b_h(\boldsymbol{p}, \boldsymbol{v}_h)
     + a_h^I((u^s,\partial_t u^b),  (\bar{v}_h^s,\bar{v}_h^b))
     + b_h^I (p^p, (\bar{v}_h^s,\bar{v}_h^b))
     = (f, v_h)_{\Omega}.     
   \end{equation}
   Next, by \cref{eq:stokesbiot_2a,eq:stokesbiot_2b},
   \begin{equation}
     \label{eq:consis_ahbh-3}
     b_h(\boldsymbol{q}_h, \boldsymbol{u}) + c_h((p^p,p^b), q_h^b)
     = - (q_h, \nabla \cdot u)_{\Omega} + (q_h^b, \lambda^{-1}(\alpha p^p-p^b))_{\Omega^b} = 0.
   \end{equation}    
   Integration by parts and using \cref{eq:stokesbiot_4},
   \begin{equation}
     \label{eq:consis_bhb}
     (\kappa^{-1}z, w_h)_{\Omega^b} + b_h^b(\boldsymbol{p}^p, (w_h,0))
     =(\kappa^{-1}z+\nabla p^p,w_h)_{\Omega^b} =0. 
   \end{equation}    
   Finally, using \cref{eq:bics-zb}, $\bar{q}_h^b=0$ on $\Gamma_P^b$,
   $n^b = -n$, \cref{eq:stokesbiot_3} and \cref{eq:bc_I_u},
   \begin{equation*}
     \begin{split}
       b_h^b(\boldsymbol{q}_h^p, (z,0))
       &= -(q_h^p,\nabla \cdot z)_{\Omega^b} + \langle \bar{q}_h^p, z\cdot n^b\rangle_{\Gamma_I}
       \\
       &= (q_h^p, c_0 \partial_t p^p + \lambda^{-1}\alpha(\alpha \partial_t p^p -\partial_t p^b)-g^b)_{\Omega^b}
       - \langle \bar{q}_h^p, (u^s-\partial_t u^b)\cdot n\rangle_{\Gamma_I},
     \end{split}
   \end{equation*}    
   which, after some rearranging can be written as
   \begin{multline}
     \label{eq:consis_bhb-2}
     (c_0\partial_t p^p, q_h^p)_{\Omega^b}
     + c_h((\partial_tp^p,\partial_tp^b), \alpha q_h^p)
     - b_h^{b}(\boldsymbol{q}_h^p, (z, 0))
      \\
     - b_h^I(\bar{q}_h^p,(u^s, \partial_tu^b))
     = (g^b, q_h^p)_{\Omega^b}.
   \end{multline}
   The result follows after comparing
   \cref{eq:consis_ahbh-2,eq:consis_ahbh-3,eq:consis_bhb,eq:consis_bhb-2}
   to \cref{eq:semidiscreteHDG}. 
\end{proof}

Before demonstrating stability of the semi-discrete problem \cref{eq:semidiscreteHDG} and
well-posedness of the fully-discrete problem
\cref{eq:fullydiscreteHDG} we first introduce some preliminary
notation and results. We start by defining:
\begin{align*}
  V^j &:= \cbr[1]{v \in [H^2(\Omega^j)]^{d}\ :\ v|_{\Gamma_D^j} = 0} \cap H^2(\Omega^j)^{d}, \quad j=s,b,
  \\
  Q^j &:= H^1(\Omega^j), \quad j=s,b, 
  \\
  Z &:= \cbr[1]{ v \in [H^1(\Omega^b)]^{d}\ :\ v \cdot n|_{\Gamma_F^b} = 0},
  \\
  Q^{b0} &:= \cbr[1]{ q \in H^2(\Omega^b) \ : \ q|_{\Gamma_P^b} = 0 } \cap H^2(\Omega^b).
\end{align*}
Let $\bar{V}^j$ and $\bar{Q}^j$ be the trace spaces of, respectively,
$V^j$ and $Q^j$ restricted to $\Gamma_0^j$ ($j=s,b$), and let
$\bar{Q}^{b0}$ be the trace space of $Q^{b0}$ restricted to
$\Gamma_0^b$. Where no confusion can occur we will write $u^j$
restricted to $\Gamma_0^j$ as $u^j$ instead of $u^j|_{\Gamma_0^j}$,
and similarly for the other unknowns.

Following the notation used also for the discrete spaces, we write for
$j=s,b$ $\boldsymbol{V}^j := V^j \times \bar{V}^j$ and
$\boldsymbol{Q}^j := Q^j \times \bar{Q}^j$. Furthermore,
$\boldsymbol{Q}^{b0} := Q^{b0} \times \bar{Q}^{b0}$. Extended function
spaces are defined as:
\begin{align*}
  \boldsymbol{V}^j(h) &:= \boldsymbol{V}_h^j + \boldsymbol{V}^j, && 
  \boldsymbol{Q}^j(h) := \boldsymbol{Q}_h^j + \boldsymbol{Q}^j, && j=s,d,
  \\
  \boldsymbol{Q}^{b0}(h) &:= \boldsymbol{Q}_h^{b0} + \boldsymbol{Q}^{b0}, &&
  Z(h) := Z_h + Z, &&
\end{align*}
with the following norms defined on $\boldsymbol{V}^j(h)$ and
$\boldsymbol{Q}^j(h)$:
\begin{align*}
  \tnorm{\boldsymbol{v}^j}_{v,j}^2&:=\sum_{K\in \mathcal{T}^j}\del[1]{\norm[0]{\varepsilon(v^j)}_K^2+h_K^{-1}\norm[0]{v^j-\bar{v}^j}_{\partial K}^2}&&\forall \boldsymbol{v}^j\in \boldsymbol{V}^j(h), && j=s,d,
  \\
  \tnorm{\boldsymbol{v}^j}_{v',j}^2&:=\tnorm{\boldsymbol{v}^j}_{v,j}^2+\sum_{K\in \mathcal{T}^j} h_K^2 |v^j|_{2,K}^2 &&\forall \boldsymbol{v}^j\in \boldsymbol{V}^j(h), && j=s,d,
  \\
  \tnorm{\boldsymbol{q}}_{q,j}^2&:=\norm[0]{q}_{\Omega^j}^2+\sum_{K \in \mathcal{T}^j}h_K\norm[0]{\bar{q}^j}_{\partial K}^2&&\forall \boldsymbol{q} \in \boldsymbol{Q}^j(h), && j=s,d.
\end{align*}
For functions $\boldsymbol{v}_h \in \boldsymbol{V}_h$ and
$\boldsymbol{q}_h \in \boldsymbol{Q}_h$ it will be useful to also
define
\begin{equation*}
  \begin{split}
    \tnorm{\boldsymbol{v}_h}_{v}^2
    &:= \tnorm{\boldsymbol{v}_h^s}_{v,s}^2+\tnorm{\boldsymbol{v}_h^b}_{v,b}^2 + \norm[0]{(\bar{v}^s_h - \bar{v}^b_h)^t}_{\Gamma_I}^2,
    \\
    \tnorm{\boldsymbol{q}_h}_q^2
    &:=\tnorm{\boldsymbol{q}_h^s}_{q,s}^2+\tnorm{\boldsymbol{q}_h^b}_{q,b}^2.        
  \end{split}
\end{equation*}

In what follows, $C > 0$ will denote a constant that is independent of
$h$ and $\Delta t$. A consequence of \cite[Lemmas 1 and
3]{Cesmelioglu:2020} are the following inequalities:
\begin{subequations}
  \begin{align}
    \label{eq:ah-cont-j}
    a_h^j(\boldsymbol{u}, \boldsymbol{v})
    &\leq C\mu^j\tnorm{\boldsymbol{u}}_{v',j}\tnorm{\boldsymbol{v}}_{v',j}
    && \forall \boldsymbol{u}, \boldsymbol{v}\in \boldsymbol{V}^j(h),
    \\
    \label{eq:ah-coer}
    a_h^j(\boldsymbol{v}_h^j, \boldsymbol{v}_h^j)
    &\geq C\mu^j\tnorm{\boldsymbol{v}_h^j}_{v,j}^2
    && \forall \boldsymbol{v}_h^j\in \boldsymbol{V}^j_h.
  \end{align}  
\end{subequations}
Due to the equivalence between $\tnorm{\cdot}_{v',j}$ and
$\tnorm{\cdot}_{v,j}$ on $\boldsymbol{V}_h^j$,
$\tnorm{\boldsymbol{u}}_{v',j}$ in \cref{eq:ah-cont-j} can be replaced
with $\tnorm{\boldsymbol{u}}_{v,j}$ if $\boldsymbol{u}$ belongs to
$\boldsymbol{V}^j_h$ (and similarly if $\boldsymbol{v}$ belongs to
$\boldsymbol{V}^j_h$). Note that, as typical of interior penalty
methods, that \cref{eq:ah-coer} holds only for a large enough
$\beta^j$.

By the Cauchy--Schwarz and Korn's inequalities, we have the
following boundedness result for $b_h^j$:
\small{ 
\begin{equation}
  \label{ineq:bhj-bound}
  \begin{split}
    b_h^j(\boldsymbol{p}^j, \boldsymbol{v}^j)
    &\leq \norm[0]{p^j}_{\Omega^j}\norm[0]{\nabla v^j}_{\Omega^j}
    + \del[1]{\sum_{K\in \mathcal{T}^j} h_K \norm[0]{\bar{p}^j}_{\partial K}^2}^{1/2}
    \del[1]{\sum_{K\in \mathcal{T}^j} h_K^{-1} \norm[0]{v^j-\bar{v}^j}_{\partial K}^2}^{1/2}
    \\
    &\leq C\del[1]{\norm[0]{p^j}_{\Omega^j}^2 + \sum_{K\in \mathcal{T}^j} h_K \norm[0]{\bar{p}^j}_{\partial K}^2}^{1/2}
    \del[1]{\norm[0]{\varepsilon(v^j)}_{\Omega^j}^2+\sum_{K\in \mathcal{T}^j} h_K^{-1} \norm[0]{v^j-\bar{v}^j}_{\partial K}^2}^{1/2}
    \\
    & \leq C\tnorm{\boldsymbol{p}^j}_{q,j} \tnorm{\boldsymbol{v}^j}_{v,j}
    \quad \forall \boldsymbol{p}^j\in \boldsymbol{Q}^j_h\quad \forall \boldsymbol{v}^j\in \boldsymbol{V}^j(h).    
  \end{split}
\end{equation}
}
Next, we discuss various inf-sup conditions on $b_h$ and $b_h^j$ that are fundamental in our proofs. First, for $j=s,b$
\begin{equation}
  \label{eq:total-pressure-inf-sup-alt}
  \inf_{\boldsymbol{0} \ne \boldsymbol{q}_h \in \boldsymbol{Q}_h^j, \boldsymbol{q}_h}
  \sup_{\boldsymbol{0} \ne \boldsymbol{v}_h \in \widetilde{\boldsymbol{V}}_h^j}
  \frac{b_h^j(\boldsymbol{q}_h, \boldsymbol{v}_h)}{\tnorm{\boldsymbol{v}_h}_{v,j} \tnorm{\boldsymbol{q}_h}_{q,j} } \ge C,
\end{equation}
where
$\widetilde{\boldsymbol{V}}^j_h := \cbr[1]{ \boldsymbol{v}_h \in
  \boldsymbol{V}_h^j\,:\, \; \bar{v}_h|_{\Gamma_I} = 0}$ is a subspace
of $\boldsymbol{V}_h^j$. The proof of
\cref{eq:total-pressure-inf-sup-alt} is given in
\cref{ap:infsupproof}. Let
$\widehat{\boldsymbol{V}}_h := \cbr[1]{ \boldsymbol{v}_h \in
  \boldsymbol{V}_h \ : \ \bar{v}_h^s \cdot n = \bar{v}_h^b \cdot n \text{ on }
  \Gamma_I}$. Then we also have
\begin{equation}
  \label{eq:inf-sup}
  \inf_{\boldsymbol{0}\neq \boldsymbol{q}_h\in \boldsymbol{Q}_h}\sup_{\boldsymbol{0}\neq \boldsymbol{v}_h\in \widehat{\boldsymbol{V}}_h}
  \frac{b_h(\boldsymbol{q}_h, \boldsymbol{v}_h)}{\tnorm{\boldsymbol{v}_h}_v\tnorm{\boldsymbol{q}_h}_q}\geq C,
\end{equation}
which was proven in \ref{ap:infsupprooffull}. In \cite[Appendix
A]{Cesmelioglu:2022} we proved
\begin{equation}
  \label{eq:infsupforwh0}
  \inf_{\boldsymbol{0} \ne \boldsymbol{q}_h^p\in \boldsymbol{Q}_h^{b,0}} \sup_{0 \ne w_h \in V_h^b}
  \frac{b_h^b(\boldsymbol{q}_h^p, (w_h,0))}{\norm[0]{w_h}_{\Omega^b} \tnorm{\boldsymbol{q}_h}_{q,b} } \ge C.
\end{equation}
We now proceed with proving stability of the semi-discrete HDG scheme \cref{eq:semidiscreteHDG}. 
\begin{theorem}[Stability]
    Let $f^b\in W^{1,1}(J;L^2(\Omega^b))$, $f^s\in L^2(J;L^2(\Omega^s))$, and $g^b\in L^2(J;L^2(\Omega^b))$ and suppose that $(\boldsymbol{u}_h, \boldsymbol{p}_h, z_h,\boldsymbol{p}_h^p) \in \mathcal{C}^1(J;\boldsymbol{X}_h)$ is a solution to \cref{eq:semidiscreteHDG}. Then 
    \small{
    \begin{align}
        \label{eq:stab-1}     
          X(t) &\leq   X(0) + C\Big[(\mu^s)^{-1/2}\|f^s (s)\|_{L^2(0,t;L^2(\Omega^s))} + \kappa^{-1/2}\|g^b(s)\|_{L^2(0,t;L^2(\Omega^b))} 
          \\
          &\hspace{7em}+(\mu^b)^{-1/2}\big(\|\partial_t f^b(s)\|_{L^1(0,t;L^2(\Omega^b))} \dif s + \max_{0\leq s\leq t}\|f^b(s)\|_{\Omega^b}\big)\Big],\nonumber
          \\
          \label{eq:stab-2}      
          \Big(\int_0^t Y(s)^2\dif s\Big)^{1/2} &\leq C\Big[X(0)+(\mu^s)^{-1/2}\|f^s\|_{L^2(0,t;L^2(\Omega^s))}+\kappa^{-1/2}\|g^b\|_{L^2(0,t;L^2(\Omega^b))} \\
          &\hspace{7em}+(\mu^b)^{-1/2}\big(\|\partial_t f^b(s)\|_{L^1(0,t;L^2(\Omega^b))} \dif s + \max_{0\leq s\leq t}\|f^b(s)\|_{\Omega^b}\big)\Big], \nonumber
    \end{align} 
    }
    where 
    \begin{equation}
        \label{eq:XY}
        \begin{split}
            X(t)^2&:= a_h^b(\boldsymbol{u}_h^b(t), \boldsymbol{u}_h^b(t)) + c_0\norm[0]{p_h^p(t)}^2_{\Omega^b}+\lambda^{-1}\norm[0]{\alpha p_h^p(t) - p_h^b(t)}_{\Omega^b}^2,
            \\
            Y(t)^2&:=a_h^s(\boldsymbol{u}_h^s(t),\boldsymbol{u}_h^s(t))
             +\gamma(\mu^s/\kappa)^{1/2}\|(\bar{u}_h^s(t)-\partial_t\bar{u}_h^b(t))^t\|_{\Gamma_I}^2 +  \kappa^{-1}\|z_h(t)\|_{\Omega^b}^2,
        \end{split}
    \end{equation}
    and $C>0$ is a constant independent of $t>0$.
\end{theorem}
\begin{proof}
    After differentiating in time the Biot part of \cref{eq:semidiscreteHDG-b}, we let $\boldsymbol{v}_h^s=\boldsymbol{u}_h^s$ , $\boldsymbol{v}_h^b=\partial_t\boldsymbol{u}_h^b$, $\boldsymbol{q}_h=-\boldsymbol{p}_h$, $\boldsymbol{q}_h^p=\boldsymbol{p}_h^p$, and $w_h=z_h$ in \cref{eq:semidiscreteHDG} and add the resulting equations. This yields
    \begin{multline}
    \label{eq:stabproof-1}
          \dfrac12 \dfrac{d}{dt}\big[a_h^b(\boldsymbol{u}_h^b, \boldsymbol{u}_h^b) + c_0\|p_h^p\|_{\Omega^b}  +\lambda^{-1}\|\alpha p_h^p-p_h^b\|_{\Omega^b}^2\textbf{}\big]  
          + \big[a_h^s(\boldsymbol{u}_h^s,\boldsymbol{u}_h^s)
          +\gamma(\mu^s/\kappa)^{1/2}\|(\bar{u}_h^s-\partial_t\bar{u}_h^b)^t\|_{\Gamma_I}^2
          \\ +  \kappa^{-1}\|z_h\|_{\Omega^b}^2  \big]
          = (f^s, u_h^s)_{\Omega^s} +  (f^b, \partial_t u_h^b)_{\Omega^b} +  (g^b, p_h^p)_{\Omega^b}.
    \end{multline}
    Observe that by the discrete Korn's inequality, coercivity of $a_h^j$, $j=s,b$ \cref{eq:ah-coer}, and the inf-sup condition \cref{eq:infsupforwh0}, 
    the following inequalities hold:
    \begin{align*}
        \|u_h^s\|_{\Omega^s}& \leq C\tnorm{\boldsymbol{u}_h^s}_{v,s}\leq C(\mu^s)^{-1/2}a_h^s(\boldsymbol{u}_h^s,\boldsymbol{u}_h^s)^{1/2}\leq C(\mu^s)^{-1/2}Y(t),
        \\
        \|u_h^b\|_{\Omega^b}& \leq C\tnorm{\boldsymbol{u}_h^b}_{v,b}\leq C(\mu^b)^{-1/2}a_h^b(\boldsymbol{u}_h^b,\boldsymbol{u}_h^b)^{1/2}\leq C(\mu^b)^{-1/2}X(t),
        \\
        \|p_h^p\|&\leq \tnorm{\boldsymbol{p}_h^p}_{q,b}\leq C\sup_{0\neq w_h\in V_h^b}\dfrac{(\kappa^{-1}z_h,w_h)_{\Omega^b}}{\|w_h\|_{\Omega^b}}\leq  C\kappa^{-1}\|z_h\|_{\Omega^b}\leq C\kappa^{-1/2}Y(t) .
    \end{align*}
    Integrating \cref{eq:stabproof-1} from 0 to $t$, $t\geq 0$, and using the above inequalites and \cref{eq:XY} in combination with the Cauchy--Schwarz inequality, we obtain:
    \begin{equation}
    \label{eq:stabproof-2}
        \begin{split}
            \dfrac12\big(X(t)^2&-X(0)^2\big) + \int_0^t Y(s)^2 \dif s 
            \\
            \leq & C(\mu^s)^{-1/2}\int_0^t \|f^s(s)\|_{\Omega^s}Y(s) \dif s + C \kappa^{-1/2}\int_0^t \|g^b(s)\|_{\Omega^b}Y(s) \dif s
            \\
            & +  C(\mu^b)^{-1/2}\big(\int_0^t \|\partial_t f^b(s)\|_{\Omega^b}X(s) \dif s + \|f^b(t)\|_{\Omega^b}X(t) + \|f^b(0)\|_{\Omega^b}X(0)\big),
        \end{split}
    \end{equation}
    where to obtain the third term on the right hand side of \cref{eq:stabproof-2} we applied integration by parts before applying the Cauchy--Schwarz inequality.
    Let us first prove \cref{eq:stab-1,eq:stab-2} under the assumption  
    \begin{equation}
    \label{eq:stab-assump}
        \max_{0\leq s\leq t} X(s) = X(t)>0.
    \end{equation}
    This assumption and Young's inequality allow us to rewrite \cref{eq:stabproof-2} as
    \begin{multline}
    \label{eq:stabproof-3}
        X(t)^2+\int_0^t Y(s)^2 \dif s 
        \\
        \leq X(0)^2 + C\alpha(t)^2 + C (\mu^b)^{-1/2}\big(\int_0^t \|\partial_t f^b(s)\|_{\Omega^b} \dif s + \|f^b(t)\|_{\Omega^b} + \|f^b(0)\|_{\Omega^b}\big)X(t)
    \end{multline}
    where $\alpha(t):=\big((\mu^s)^{-1}\int_0^t \|f^s(s)\|^2_{\Omega^s}\dif s + \kappa^{-1}\int_0^t \|g^b(s)\|^2_{\Omega^b}\dif s\big)^{1/2}$. Here we remark that the constant in Young's inequality is independent of $t$.
    
    If $\alpha(t)>X(t)$, inequality \cref{eq:stab-1} is obvious. Assume therefore that $\alpha(t)\leq X(t)$. Recalling that $X(0)\leq X(t)$ by \cref{eq:stab-assump}, and dividing both sides of \cref{eq:stabproof-3} by $X(t)>0$, we obtain 
    \begin{equation}
    \label{eq:XtboundX0etc}
            X(t) \leq X(0) + C \alpha(t) +  C(\mu^b)^{-1/2}\big(\int_0^t \|\partial_t f^b(s)\|_{\Omega^b} \dif s + \|f^b(t)\|_{\Omega^b} + \|f^b(0)\|_{\Omega^b}\big),
    \end{equation}
    implying \cref{eq:stab-1}. \Cref{eq:stab-2} follows from \cref{eq:stabproof-3,eq:stab-1}.
        
    If assumption \cref{eq:stab-assump} is not true, there exists $0 \le \tilde{t} < t$ such that
    \begin{equation}
    \label{eq:stab-assump-nn}
        \max_{0 \le s \le t} X(s) = X(\tilde{t}) > 0.
    \end{equation}
    Integrating \cref{eq:stabproof-1} from 0 to $\tilde{t}$, using the same steps as above that were used to find \cref{eq:stab-1}, but restricting ourselves to the time interval $(0,\tilde{t})$, and using \cref{eq:stab-assump-nn} instead of \cref{eq:stab-assump}, we find
    \begin{multline}
    \label{eq:stab-1-nn}
        X(\tilde{t}) \leq   X(0) + C\Big[(\mu^s)^{-1/2}\|f^s (s)\|_{L^2(0,\tilde{t};L^2(\Omega^s))} + \kappa^{-1/2}\|g^b(s)\|^2_{L^2(0,\tilde{t};L^2(\Omega^b))} 
        \\
        +(\mu^b)^{-1/2}\big(\|\partial_t f^b(s)\|_{L^1(0,\tilde{t};L^2(\Omega^b))} \dif s + \max_{0\leq s\leq \tilde{t}}\|f^b(s)\|_{\Omega^b}\big)\Big].
    \end{multline}    
    Then, \cref{eq:stab-1} holds because $X(t) < X(\tilde{t})$ by assumption \cref{eq:stab-assump-nn} and $\tilde{t} < t$. Finally, to prove \cref{eq:stab-2}, note that \cref{eq:stabproof-2} holds for any $t$, so a crude inequality by Young's inequality and the assumption \cref{eq:stab-assump-nn} give 
    \begin{align*}
        \int_0^t Y(s)^2 \dif s &\leq X(0)^2 + C\big((\mu^s)^{-1}\int_0^t \|f^s(s)\|^2_{\Omega^s}\dif s + \kappa^{-1}\int_0^t \|g^b(s)\|^2_{\Omega^b}\dif s\big) \\
        &\quad + C (\mu^b)^{-1/2}\big(\int_0^t \|\partial_t f^b(s)\|_{\Omega^b} \dif s + \max_{0\le s\le t} \|f^b(s)\|_{\Omega^b} \big)X(\tilde{t}) \\
        &\leq X(0)^2 + C\big((\mu^s)^{-1}\int_0^t \|f^s(s)\|^2_{\Omega^s}\dif s + \kappa^{-1}\int_0^t \|g^b(s)\|^2_{\Omega^b}\dif s\big) \\
        &\quad + C^2 (\mu^b)^{-1}\big(\int_0^t \|\partial_t f^b(s)\|_{\Omega^b} \dif s + \max_{0\le s\le t} \|f^b(s)\|_{\Omega^b} \big)^2 + \frac 14 X(\tilde{t})^2 . 
    \end{align*}
    \Cref{eq:stab-2} now follows by combining this result with \cref{eq:stab-1-nn} and using that $\tilde{t} < t$. 
\end{proof}

Well-posedness of the fully discrete method is now given by the
following lemma.
\begin{lemma}[Existence and uniqueness]
  \label{lem:existenceuniqueness}
  There exists a unique solution to the fully-discrete HDG method
  \cref{eq:fullydiscreteHDG}.
\end{lemma}
\begin{proof}
  Setting the right hand sides of \cref{eq:fullydiscreteHDG} to zero,
  $(\boldsymbol{u}_h^{n}, \boldsymbol{p}_h^{n}, z_h^{n},
  \boldsymbol{p}_h^{p,n})=(\boldsymbol{0}, \boldsymbol{0}, 0,
  \boldsymbol{0})$, and adding the resulting equations, we obtain
  \begin{equation}
    \label{eq:zerodataeqn}
    \begin{split}
      &a_h(\boldsymbol{u}_h^{n+1}, \boldsymbol{v}_h) 
      + a_h^I((\bar{u}_h^{s,n+1}, \tfrac{1}{\Delta t}\bar{u}_h^{b,n+1}),(\bar{v}_h^s,\bar{v}_h^b))
      + b_h(\boldsymbol{p}_h^{n+1}, \boldsymbol{v}_h)
      \\
      &+ b_h^I(\bar{p}_h^{p,n+1},(\bar{v}_h^s,\bar{v}_h^b))
      + b_h(\boldsymbol{q}_h, \boldsymbol{u}_h^{n+1})
      + c_h((p_h^{p,n+1},p_h^{b,n+1}),q_h^b)
      \\
      &+ \tfrac{1}{\Delta t}(c_0 p_h^{p,n+1}, q_h^p)_{\Omega^b}
      + \tfrac{1}{\Delta t}c_h((p_h^{p,n+1}, p_h^{b,n+1}), \alpha q_h^p)
      - b_h^{b}(\boldsymbol{q}_h^p, (z_h^{n+1}, 0))
      \\ 
      &- b_h^I(\bar{q}_h^p,(\bar{u}_h^{s,n+1}, \tfrac{1}{\Delta t}\bar{u}_h^{b,n+1}))
      + (\kappa^{-1} z_h^{n+1}, w_h)_{\Omega^b} + b_h^{b}(\boldsymbol{p}_h^{p,n+1}, (w_h, 0))
      = 0,        
    \end{split}
  \end{equation}
  which holds for all
  $(\boldsymbol{v}_h, \boldsymbol{q}_h, w_h, \boldsymbol{q}_h^p)\in
  \boldsymbol{X}_h$. Choosing now
  $\boldsymbol{v}_h^s = \boldsymbol{u}_h^{s,n+1}$,
  $\boldsymbol{v}_h^b = \tfrac{1}{\Delta t}\boldsymbol{u}_h^{b,n+1}$,
  $\boldsymbol{q}_h^s = -\boldsymbol{p}_h^{s,n+1}$,
  $\boldsymbol{q}_h^b = -\tfrac{1}{\Delta t}\boldsymbol{p}_h^{b,n+1}$,
  $w_h=z_h^{n+1}$, and
  $\boldsymbol{q}_h^p = \boldsymbol{p}_h^{p,n+1}$, we find
  \begin{equation*}
    \begin{split}
      &a_h^s(\boldsymbol{u}_h^{s,n+1}, \boldsymbol{u}_h^{s,n+1})
      + \tfrac{1}{\Delta t}a_h^b(\boldsymbol{u}_h^{b,n+1}, \boldsymbol{u}_h^{b,n+1})
      + \gamma(\mu^s/\kappa)^{1/2}\norm[1]{(\bar{u}_h^{s,n+1}-\frac{1}{\Delta t}\bar{u}_h^{b,n+1})^t}_{\Gamma_I}^2
      \\
      &+ \tfrac{c_0}{\Delta t}\norm[0]{p_h^{p,n+1}}_{\Omega^b}^2
      + \tfrac{\lambda^{-1}}{\Delta t}\norm[0]{\alpha p_h^{p,n+1}-p_h^{b,n+1}}_{\Omega^b}^2
      + \kappa^{-1}\norm[0]{z_h^{n+1}}_{\Omega^b}^2
      = 0.      
    \end{split}
  \end{equation*}
  Coercivity of $a_h^j$ \cref{eq:ah-coer} and positivity of $\gamma$,
  $\mu^s$, $\kappa$, $\lambda$, $\alpha$, and $\Delta t$ and
  nonnegativity of $c_0$ imply that $\boldsymbol{u}_h^{n+1}$ and
  $z_h^{n+1}$ are zero.  Substituting now $z_h^{n+1}=0$ and setting
  $\boldsymbol{v}_h=\boldsymbol{0}$,
  $\boldsymbol{q}_h=\boldsymbol{0}$, and
  $\boldsymbol{q}_h^p=\boldsymbol{0}$ in \cref{eq:zerodataeqn}, we
  obtain $b_h^{b}(\boldsymbol{p}_h^{p,n+1}, (w_h, 0)) = 0$ for all
  $w_h \in V_h^b$. It follows from \cref{eq:infsupforwh0} that
  $\boldsymbol{p}_h^{p,n+1} = \boldsymbol{0}$. Next, substituting
  $\boldsymbol{u}_h^{n+1} = \boldsymbol{0}$, choosing
  $\bar{v}_h^s=\bar{v}_h^b=0$ on $\Gamma_I$, and choosing $w_h=0$,
  $\boldsymbol{q}_h=\boldsymbol{0}$,
  $\boldsymbol{q}_h^p=\boldsymbol{0}$ in \cref{eq:zerodataeqn}, we
  obtain $b_h^s(\boldsymbol{p}_h^{s,n+1}, \boldsymbol{v}_h^s)=0$ for
  all $\boldsymbol{v}_h^s \in \widetilde{\boldsymbol{V}}_h^s$ and
  $b_h^b(\boldsymbol{p}_h^{b,n+1}, \boldsymbol{v}_h^b)=0$ for all
  $\boldsymbol{v}_h^b \in \widetilde{\boldsymbol{V}}_h^b$. It follows
  now from \cref{eq:total-pressure-inf-sup-alt} that
  $\boldsymbol{p}_h^{j,n+1}=\boldsymbol{0}$ for $j=s,b$. 
\end{proof}

\section{A priori error analysis}
\label{sec:apriorierrorhdg}

We now prove convergence results for the HDG method. Throughout this
section the superscript $j$ will refer to $s$ and $b$. For this we
first introduce suitable interpolation/projection operators. For
vector valued functions we define $\Pi_V:[H^1(\Omega^b)]^d \to V_h^b$
to be the Brezzi--Douglas--Marini (BDM) interpolation operator
\cite[Section III.3]{Brezzi:book} and
$\boldsymbol{\Pi}_V^{\text{ell},j} := (\Pi_V^{\text{ell},j},
\bar{\Pi}_V^{\text{ell},j}) : [H^1(\Omega^j)]^d \to
\boldsymbol{V}_h^j$ to be the elliptic interpolation operator defined
as
\begin{equation*}
  a_h^j(\boldsymbol{\Pi}_{V}^{\text{ell},j}u, \boldsymbol{v}_h) 
  =a_h^j((u,u), \boldsymbol{v}_h), \qquad \forall \boldsymbol{v}_h \in \boldsymbol{V}_h^j.
\end{equation*}
For scalar functions we denote by $\Pi_Q^j$, $\bar{\Pi}_Q^j$, and
$\bar{\Pi}_{Q^{0}}^b$ the $L^2$-projection operators onto,
respectively, $Q_h^j$, $\bar{Q}_h^j$, and $\bar{Q}_h^{b0}$. We next
introduce the following notation for the errors:
\begin{subequations}
  \label{eq:errorsplitting}
  \begin{align}
    u^j - u_h^j &= (u^j - \Pi_V^{\text{ell},j}u^j) - (u_h^j - \Pi_V^{\text{ell},j}u^j) &&= e_{u^j}^{I} - e_{u^j}^{h},
    \\
    u^j|_{\Gamma_0^j} - \bar{u}_h^j &= (u^j|_{\Gamma_0^j} - \bar{\Pi}_V^{\text{ell},j}u^j) - (\bar{u}_h^j - \bar{\Pi}_V^{\text{ell},j}u^j) &&= \bar{e}_{u^j}^{I} - \bar{e}_{u^j}^{h},
    \\
    z - z_h &= (z - \Pi_Vz) - (z_h - \Pi_Vz) &&= e_z^{I} - e_z^{h},
    \\
    p^j - p_h^j &= (p^j - \Pi_Q^jp^j) - (p_h^j - \Pi_Q^jp^j) &&= e_{p^j}^{I} - e_{p^j}^{h},
    \\
    p^j|_{\Gamma_0^j} - \bar{p}_h^j &= (p^j|_{\Gamma_0^j} - \bar{\Pi}_Q^jp^j) - (\bar{p}_h^j - \bar{\Pi}_Q^jp^j) &&= \bar{e}_{p^j}^{I} - \bar{e}_{p^j}^{h},
    \\
    p^p - p_h^p &= (p^p - \Pi_Q^bp^p) - (p_h^p - \Pi_Q^bp^p) &&= e_{p^p}^{I} - e_{p^p}^{h},
    \\
    p^p|_{\Gamma_0^b} - \bar{p}_h^p &= (p^p|_{\Gamma_0^b} - \bar{\Pi}_{Q^0}^bp^p) - (\bar{p}_h^p - \bar{\Pi}_{Q^0}^bp^p) &&= \bar{e}_{p^p}^{I} - \bar{e}_{p^p}^{h}.
  \end{align}  
\end{subequations}
We also define $e_u^I$ and $e_p^I$ such that
$e_u^I|_{\Omega^j} = e_{u^j}^I$ and $e_p^I|_{\Omega^j} =
e_{p^j}^I$. Similarly, $e_u^h$ and $e_p^h$ are defined such that
$e_u^h|_{\Omega^j} = e_{u^j}^h$ and $e_p^h|_{\Omega^j} =
e_{p^j}^h$. 

We then write
\begin{equation*}
  \boldsymbol{e}_u^I
  = (e_u^I, \bar{e}_{u^s}^I, \bar{e}_{u^b}^I),   
  \quad
  \boldsymbol{e}_u^h = (e_u^h, \bar{e}_{u^s}^h, \bar{e}_{u^b}^h),
  \quad
  \boldsymbol{e}_p^I
  = (e_p^I, \bar{e}_{p^s}^I, \bar{e}_{p^b}^I),
  \quad
  \boldsymbol{e}_p^h = (e_p^h, \bar{e}_{p^s}^h, \bar{e}_{p^b}^h).  
\end{equation*}
It will furthermore be convenient to introduce the following notation:
\begin{align*}
  \boldsymbol{e}_{u^j}^{I} &= (e_{u^j}^{I}, \bar{e}_{u^j}^{I}),
  && \boldsymbol{e}_{p^j}^{I} = (e_{p^j}^{I}, \bar{e}_{p^j}^{I}),
  && \boldsymbol{e}_{p^p}^{I} = (e_{p^p}^{I}, \bar{e}_{p^p}^{I}),
  \\
  \boldsymbol{e}_{u^j}^{h} &= (e_{u^j}^{h}, \bar{e}_{u^j}^{h}),
  && \boldsymbol{e}_{p^j}^{h} = (e_{p^j}^{h}, \bar{e}_{p^j}^{h}),
  && \boldsymbol{e}_{p^p}^{h} = (e_{p^p}^{h}, \bar{e}_{p^p}^{h}).
\end{align*}
The following interpolation estimates hold:
\begin{subequations}
  \begin{align}
    \label{eq:interpolation-property-bdm}
    \norm[0]{e_z^I}_{m,K}
    &\leq Ch_K^{\ell-m} \norm[0]{z}_{\ell,K},
      \quad m=0, 1, 2, \quad \max\{m,1\} \leq \ell \leq k+1,
    \\
    \label{eq:interpolation-property-ell}
    a_h^j( \boldsymbol{e}_{u^j}^{I}, \boldsymbol{e}_{u^j}^{I} )^{\frac 12}
    &\leq C \sqrt{\mu^j} h^{\ell-1}\norm[0]{u^j}_{\ell,\Omega^j},
      \quad 1\leq \ell\leq k+1, \quad j=s,b,
  \end{align}  
\end{subequations}
where \cref{eq:interpolation-property-bdm} is the usual BDM
interpolation estimate \cite[Lemma 7]{Hansbo:2002} and
\cref{eq:interpolation-property-ell} follows from standard a priori
error estimate theory for second order elliptic equations.

To initialize \cref{eq:fullydiscreteHDG}, we set
\begin{equation}
  \label{eq:initialconditions0}
    u_h^{b,0} = \Pi_V^{\text{ell},j}u^b(0), \quad
    \bar{u}_h^{b,0} = \bar{\Pi}_V^{\text{ell},b}u^b(0), \quad
    p_h^{p,0} = \Pi_Q^bp^p(0).
\end{equation}
The initial total pressure is set by
$p_h^{b,0} = \alpha p_h^{p,0} - \lambda \nabla \cdot u_h^{b,0}$.

\begin{theorem}[Error equation]
  \label{thm:erroreqn}
  Let
  $(\boldsymbol{u}_h^{n}, \boldsymbol{p}_h^{n}, z_h^{n},
  \boldsymbol{p}_h^{p,n})$, for $n=1,\hdots,N$, be the solution to
  \cref{eq:fullydiscreteHDG} with initial conditions
  \cref{eq:initialconditions0}. Let $(u,p,z,p^p)$ be a solution to the
  coupled Stokes--Biot problem
  \cref{eq:stokesbiot,eq:interface,eq:bics} on time interval $J=(0,T]$
  and let $\boldsymbol{u} := (u, u|_{\Gamma_0^s}, u|_{\Gamma_0^b})$,
  $\boldsymbol{p} := (p, p|_{\Gamma_0^s}, p|_{\Gamma_0^b})$, and
  $\boldsymbol{p}^p := (p^p, p^p|_{\Gamma_0^b})$. Then, with the exact
  solution evaluated at $t=t^{n+1}$ (and the superscript $n+1$
  suppressed for notational convenience),
  \begin{subequations}
    \label{eq:erroreqn}
    \begin{align}
      \label{eq:erroreqn-a}
      &a_h(\boldsymbol{e}_u^h, \boldsymbol{v}_h) 
        + a_h^I((\bar{e}_{u^s}^h,d_t\bar{e}_{u^b}^h),(\bar{v}_h^s,\bar{v}_h^b))
        + b_h(\boldsymbol{e}_p^h, \boldsymbol{v}_h)
        + b_h^I(\bar{e}_{p^p}^h,(\bar{v}_h^s,\bar{v}_h^b))
      \\
      \nonumber
      &\hspace{1em}= a_h^I((0,\partial_tu^b-d_tu^b),(\bar{v}_h^s,\bar{v}_h^b))
        + a_h^I((\bar{e}_{u^s}^I,d_t\bar{e}_{u^b}^I),(\bar{v}_h^s,\bar{v}_h^b)),   
      \\
      \label{eq:erroreqn-b}
      &(\kappa^{-1} e_z^h, w_h)_{\Omega^b} + b_h^{b}(\boldsymbol{e}_{p^p}^h, (w_h, 0))
        =
        (\kappa^{-1} e_z^I, w_h)_{\Omega^b},
      \\
      \label{eq:erroreqn-c}
      &(c_0 d_te_{p^p}^h, q_h^p)_{\Omega^b} 
        + c_h((d_te_{p^p}^h, d_te_{p^b}^h), \alpha q_h^p)
      - b_h^{b}(\boldsymbol{q}_h^p, (e_z^h, 0)) 
        - b_h^I(\bar{q}_h^p,(\bar{e}_{u^s}^h, d_t\bar{e}_{u^b}^h))
      \\
      \nonumber
      & \hspace{1em} = (c_0 (\partial_t p^p - d_t p^p), q_h^p)_{\Omega^b}
        + c_h((\partial_t p^p - d_t p^p, \partial_t p^b-d_t p^b), \alpha q_h^p)
      - b_h^I(\bar{q}_h^p,(0, \partial_t u^b - d_t u^b)),
      \\
      \label{eq:erroreqn-d}
      &b_h^b(\boldsymbol{q}_h^b, d_t\boldsymbol{e}_{u^b}^h)
        + c_h((d_te_{p^p}^h,d_te_{p^b}^b), q_h^b)
       = b_h^b(\boldsymbol{q}_h^b, d_t\boldsymbol{e}_{u^b}^I),          
      \\
      \label{eq:erroreqn-e}
      &b_h^s(\boldsymbol{q}_h^s, \boldsymbol{e}_{u^s}^h)
        = b_h^s(\boldsymbol{q}_h^s, \boldsymbol{e}_{u^s}^I).
    \end{align}  
  \end{subequations}  
\end{theorem}
\begin{proof}
  By \cref{eq:fullydiscreteHDG} and \cref{lem:consistency},
  \begin{subequations}
    \label{eq:efullydiscreteHDG}
    \begin{align}
      \label{eq:efullydiscreteHDG-a}
      &a_h(\boldsymbol{u}_h-\boldsymbol{u}, \boldsymbol{v}_h) 
        + a_h^I((\bar{u}_h^{s}-u^s,d_t\bar{u}_h^{b}-\partial_tu^b),(\bar{v}_h^s,\bar{v}_h^b))
      \\ \nonumber
      &\hspace{7em} + b_h(\boldsymbol{p}_h-\boldsymbol{p}, \boldsymbol{v}_h)
        + b_h^I(\bar{p}_h^{p}-p^p,(\bar{v}_h^s,\bar{v}_h^b))
        = 0,
      \\
      \label{eq:efullydiscreteHDG-b}
      &b_h(\boldsymbol{q}_h, \boldsymbol{u}_h-\boldsymbol{u})
        + c_h((p_h^{p}-p^p,p_h^{b}-p^b),q_h^b)
        = 0,
      \\
      \label{eq:efullydiscreteHDG-c}
      &(c_0 d_t p_h^{p} - c_0 \partial_tp^p, q_h^p)_{\Omega^b}
        + c_h((d_tp_h^{p}-\partial_tp^p,d_tp_h^{b}-\partial_tp^b), \alpha q_h^p)        
      \\ \nonumber 
      & \hspace{7em}- b_h^{b}(\boldsymbol{q}_h^p, (z_h-z, 0))
        - b_h^I(\bar{q}_h^p,(\bar{u}_h^{s}-u^s, d_t\bar{u}_h^{b}-\partial_tu^b))
        = 0,
      \\
      \label{eq:efullydiscreteHDG-d}
      &(\kappa^{-1} (z_h-z), w_h)_{\Omega^b} + b_h^{b}(\boldsymbol{p}_h^{p}-\boldsymbol{p}^p, (w_h, 0))
        = 0.
    \end{align}  
  \end{subequations}
  We split \cref{eq:efullydiscreteHDG-b} into its Stokes and Biot parts as:
  \begin{subequations}
    \label{eq:efullydiscreteHDG-b-split}
    \begin{align}
      \label{eq:efullydiscreteHDG-b-i}
      b_h^b(\boldsymbol{q}_h^b, d_t(\boldsymbol{u}_h^{b}-\boldsymbol{u}^b))
      + c_h((d_tp_h^{p}-d_tp^p,d_tp_h^{b}-d_tp^b),q_h^b)
      &= 0,
      \\
      \label{eq:efullydiscreteHDG-b-ii}
      b_h^s(\boldsymbol{q}_h^s, \boldsymbol{u}_h^{s}-\boldsymbol{u}^s)
      &= 0,      
    \end{align}
  \end{subequations}
  where we applied $d_t$ to the first equation.  Noting that
  $\partial_t u^b - d_t u^b=\partial_t u^b - d_t
  \bar{\Pi}_V^{\text{ell},b}u^b- d_t\bar{e}_{u^b}^I$,
  $\partial_t p^b-d_t p^b = \partial_t p^b-d_t \Pi_Q^b p^b-
  d_te_{p^b}^I$, and
  $\partial_t p^p - d_t p^p = \partial_t p^p - d_t \Pi_Q^b p^p-
  d_te_{p^p}^I$, combining \cref{eq:efullydiscreteHDG} and
  \cref{eq:efullydiscreteHDG-b-split}, and using the error splitting
  according to \cref{eq:errorsplitting} yields
  \begin{subequations}
    \label{eq:eefullydiscreteHDG}
    \begin{align}
      \label{eq:eefullydiscreteHDG-a}
      &a_h(\boldsymbol{e}_u^h, \boldsymbol{v}_h) 
        + a_h^I((\bar{e}_{u^s}^h,d_t\bar{e}_{u^b}^h),(\bar{v}_h^s,\bar{v}_h^b))
        + b_h(\boldsymbol{e}_p^h, \boldsymbol{v}_h)
        + b_h^I(\bar{e}_{p^p}^h,(\bar{v}_h^s,\bar{v}_h^b))
      \\
      \nonumber
      &\hspace{1em}=
        a_h(\boldsymbol{e}_u^I, \boldsymbol{v}_h)
        + b_h(\boldsymbol{e}_p^I, \boldsymbol{v}_h)
        + a_h^I((\bar{e}_{u^s}^I,\partial_tu^b-d_t\bar{\Pi}_V^{\text{ell},b}u^b),(\bar{v}_h^s,\bar{v}_h^b))
      + b_h^I(\bar{e}_{p^p}^I,(\bar{v}_h^s,\bar{v}_h^b)),
      \\
      \label{eq:eefullydiscreteHDG-d}
      &(\kappa^{-1} e_z^h, w_h)_{\Omega^b} + b_h^{b}(\boldsymbol{e}_{p^p}^h, (w_h, 0))
        =
        (\kappa^{-1} e_z^I, w_h)_{\Omega^b} + b_h^{b}(\boldsymbol{e}_{p^p}^I, (w_h, 0)),
      \\
      \label{eq:eefullydiscreteHDG-c}
      &(c_0 d_te_{p^p}^h, q_h^p)_{\Omega^b} 
        + c_h((d_te_{p^p}^h, d_te_{p^b}^h), \alpha q_h^p)
      - b_h^{b}(\boldsymbol{q}_h^p, (e_z^h, 0)) 
        - b_h^I(\bar{q}_h^p,(\bar{e}_{u^s}^h, d_t\bar{e}_{u^b}^h))
      \\
      \nonumber
      &\hspace{1em} 
      =
        (c_0 (\partial_tp^p - d_t\Pi_{Q}^b p^{p}), q_h^p)_{\Omega^b}
        + c_h((\partial_tp^p-d_t\Pi_Q^bp^p, \partial_tp^b-d_t\Pi_Q^bp_h^b), \alpha q_h^p)
      \\
      \nonumber
      &\hspace{2em}- b_h^{b}(\boldsymbol{q}_h^p, (e_z^I, 0))
        - b_h^I(\bar{q}_h^p,(\bar{e}_{u^s}^I, \partial_tu^b-d_t\bar{\Pi}_V^{\text{ell,b}}u^b)),
      \\
      \label{eq:eefullydiscreteHDG-b-i}
      &b_h^b(\boldsymbol{q}_h^b, d_t\boldsymbol{e}_{u^b}^h)
        + c_h((d_te_{p^p}^h,d_te_{p^b}^b), q_h^b)
      =
         b_h^b(\boldsymbol{q}_h^b, d_t\boldsymbol{e}_{u^b}^I)
         + c_h((d_te_{p^p}^I,d_te_{p^b}^I),q_h^b),          
      \\
      \label{eq:eefullydiscreteHDG-b-ii}
      &b_h^s(\boldsymbol{q}_h^s, \boldsymbol{e}_{u^s}^h)
        = b_h^s(\boldsymbol{q}_h^s, \boldsymbol{e}_{u^s}^I).
    \end{align}  
  \end{subequations}
  Then, by definition of the chosen
  interpolation/projection operators, the terms
  $a_h(\boldsymbol{e}_u^I, \boldsymbol{v}_h)$,
  $b_h(\boldsymbol{e}_p^I, \boldsymbol{v}_h)$,
  $b_h^I(\bar{e}_{p^p}^I,(\bar{v}_h^s,\bar{v}_h^b))$,
  $b_h^{b}(\boldsymbol{e}_{p^p}^I, (w_h, 0))$,
  $b_h^{b}(\boldsymbol{q}_h^p, (e_z^I, 0))$,
  $c_h((d_te_{p^p}^I,d_te_{p^b}^I),q_h^b)$,
  $(c_0 d_te_{p^p}^I, q_h^p)_{\Omega^b}$,
  $c_h((d_te_{p^p}^I, d_te_{p^b}^I), \alpha q_h^p)$, and
  $b_h^I(\bar{q}_h^p,(\bar{e}_{u^s}^I, d_t\bar{e}_{u^b}^I))$
  vanish. The result follows. 
\end{proof}

The following lemma will be used in the proof of the a priori error
estimates of \cref{thm:apriorierrorestimate}.
\begin{lemma}
  \label{lem:pPGammaI_z}
  Let $\bar{e}_{p^p}^h$, $z$, and $z^h$ be defined as in
  \cref{thm:erroreqn}. There exists a constant $C>0$ such that:
  \begin{equation}
    \label{eq:pPGammaI_z}
    \norm[0]{\bar{e}_{p^p}^h}_{\Gamma_I} \le C \kappa^{-1} \norm[0]{z - z_h}_{\Omega^b}.
  \end{equation}
\end{lemma}
\begin{proof}
  We start by defining
  \begin{equation*}
    \norm[0]{q_h}_{1,h}^2 := \sum_{K \in \mathcal{T}^b} \norm[0]{\nabla q_h}_K^2
    + \sum_{F \in \mathcal{F}^b \setminus (\mathcal{F}_I \cup \mathcal{F}_F^b)} h_F^{-1} \norm[0]{\jump{q_h}}_F^2,
    \quad \forall q_h \in Q_h^b.
  \end{equation*}
  Since $|\Gamma_P^b| > 0$, $\norm[0]{q_h}_{1,h}$ is a norm on
  $Q_h^b$. From the trace inequality of broken functions
  \cite[Theorem~4.4]{Buffa:2009} and discrete Poincar\'{e} inequality
  \cite{Brenner:2003} we have that for all $q_h \in Q_h^b$,
  \begin{equation}
    \label{eq:boundqhgammaI}
    \norm[0]{ q_h }_{\Gamma_I}
    \le C \del[1]{\norm[0]{ q_h }_{L^1(\Omega^b)} + \norm[0]{ q_h }_{1,h}}
    \le C \norm[0]{ q^h }_{1,h}.
  \end{equation}
  Let us now consider \cref{eq:erroreqn-b}. We find that for all
  $w_h \in V_h^b$
  \begin{equation}
    \label{eq:whzzh}
    \begin{split}
      0
      &=(\kappa^{-1} (z-z_h), w_h)_{\Omega^b} - b_h^{b}(\boldsymbol{e}_{p^p}^h, (w_h, 0))
      \\
      &=(\kappa^{-1} (z-z_h), w_h)_{\Omega^b} + (e_{p^p}^h, \nabla \cdot w_h)_{\Omega^b}
      - \langle \bar{e}_{p^p}^h, w_h \cdot n^b \rangle_{\partial \mathcal{T}^b \setminus \Gamma_P^b},
    \end{split}
  \end{equation}
  where the second equality holds because $\bar{e}_{p^p}^h = 0$ on
  $\Gamma_P^b$ (since $p^p = 0$, $\bar{p}_h^p = 0$, and
  $\bar{\Pi}_{Q^0}^bp^p = 0$ on $\Gamma_P^b$).

  \textbf{Step 1.} Choose $w_h \in V_h^b \cap H(\text{div};\Omega^b)$
  such that $w_h \cdot n = \bar{e}_{p^p}^h$ on $F \in \mathcal{F}_I$
  and such that the remaining moments for the degrees of freedom of a
  BDM element vanish. By a standard scaling argument we observe that
  \begin{equation*}
    \norm[0]{w_h}_{\Omega^b} + \del[1]{ \sum_{F \in \mathcal{F}^b}h_F\norm{w_h \cdot n}^2_F }^{1/2}
    \le C h^{1/2} \norm[0]{\bar{e}_{p^p}^h}_{\Gamma_I}.
  \end{equation*}
  Substituting the above defined $w_h$ in \cref{eq:whzzh}, integrating
  by parts and using the Cauchy--Schwarz inequality and
  \cref{eq:boundqhgammaI}:
  \begin{equation}
    \label{epphIzzh}
    \begin{split}
      \norm[0]{\bar{e}_{p^p}^h}_{\Gamma_I}^2
      &=(\kappa^{-1} (z-z_h), w_h)_{\Omega^b} + \langle e_{p^p}^h, \bar{e}_{p^p}^h \rangle_{\Gamma_I}
      \le C\del[1]{\kappa^{-1}\norm[0]{z-z_h}_{\Omega^b} + \norm[0]{e_{p^p}^h}_{1,h}} \norm[0]{\bar{e}_{p^p}^h}_{\Gamma_I}.      
    \end{split}
  \end{equation}
  
  \textbf{Step 2.} By (the proof of) \cite[Lemma~2.1]{Lovadina:2006},
  there exists a $w_h^0 \in V_h^b \cap H(\text{div};\Omega^b)$ such
  that
  \begin{align*}
    \langle w_h^0 \cdot n, r \rangle_F
    &= h_F^{-1} \langle\jump{e_{p^p}^h}, r \rangle_F
    && \forall r \in P_{k-1}(F),  \quad \forall F \in \mathcal{F}^b \setminus (\mathcal{F}_I\cup \mathcal{F}_F^b),
    \\
    (w_h^0, r)_T
    &= -(\nabla e_{p^p}^h, r)_{T}
    && \forall r \in [P_{k-2}(T)]^2, \quad \forall T \in \mathcal{T}^b,
    \\
    w_h^0 \cdot n
    &= 0 \text{ on } F
    && \forall F \in \mathcal{F}_I\cup \mathcal{F}_F^b.
    \end{align*}
    Additionaly, $w_h^0$ satisfies
    \begin{equation*}
      (\nabla \cdot w_h^0, e_{p^p}^h)_{\Omega^b} = \norm[0]{ e_{p^p}^h }_{1,h}^2,
      \qquad \norm[0]{ w_h^0 }_{\Omega^b} \le C \norm[0]{ e_{p^p}^h }_{1,h}.
    \end{equation*}
    Substituting $w_h^0$ in \cref{eq:whzzh}, we obtain:
    \begin{equation}
      \label{eq:epp1hboundzz}
      \norm[0]{ e_{p^p}^h }_{1,h}^2
      = -(\kappa^{-1} (z-z_h), w_h^0)_{\Omega^b}
      \le C \kappa^{-1}\norm[0]{z - z_h}_{\Omega^b} \norm[0]{e_{p^p}^h}_{1,h}.
    \end{equation}
    The result follows after combining \cref{epphIzzh} and
    \cref{eq:epp1hboundzz}. 
\end{proof}

The following theorem provides an a priori error estimate.

\begin{theorem}[A priori error estimate]
  \label{thm:apriorierrorestimate}
  Let $(u,p,z,p^p)$ be the solution to the coupled Stokes--Biot problem
  \cref{eq:stokesbiot,eq:interface,eq:bics} on the time interval
  $J=(0,T]$ such that $u^s \in C^0(J;H^{\ell+1}(\Omega^s))$,
  $\partial_tu^b \in L^2(J;H^{\ell+1}(\Omega^b))$,
  $\partial_{tt}u^b \in L^2(J;H^1(\Omega^b))$, $p^s \in C^0(J;L^2(\Omega^s))$,
  $p^p, p^b \in W^{2,1}(J;L^2(\Omega^b))$, and
  $z \in C^0(J;H^{\ell}(\Omega^b))$. Let
  $(\boldsymbol{u}_h^{n}, \boldsymbol{p}_h^{n}, z_h^{n},
  \boldsymbol{p}_h^{p,n}) \in \boldsymbol{X}_h$, for $n=1,\hdots,N$,
  be the solution to \cref{eq:fullydiscreteHDG} with initial
  conditions \cref{eq:initialconditions0}. The following hold:
  \begin{subequations}
    \begin{align}
      \label{ineq:apriori-1}
      &(\mu^b)^{1/2}\tnorm{\boldsymbol{e}_{u^b}^{h,n}}_{v,b}
      + c_0^{1/2} \norm[0]{e_{p^p}^{h,n}}_{\Omega^b}
      +\lambda^{-1/2}\norm[0]{\alpha  e_{p^p}^{h,n} - e_{p^b}^{h,n}}_{\Omega^b} \\ \nonumber
      &+(\mu^s)^{1/2}\big(\Delta t \sum_{i=1}^n \tnorm{\boldsymbol{e}_{u^s}^{h,i}}_{v,s}^2\big)^{1/2} 
      + \kappa^{-1/2} \big(\Delta t\sum_{i=1}^n\norm[0]{e_{z}^{h,i}}^2_{\Omega^b}\big)^{1/2} \\ \nonumber
      &+ \gamma^{1/2}(\mu^s/\kappa)^{1/4} \big( \Delta t \sum_{i=1}^n \norm[0]{(\bar{e}_{u^s}^{h,i} - d_t\bar{e}_{u^b}^{h,i})^t}_{\Gamma_I}^2\big)^{1/2}
      \leq \bar{C}_1\Delta t + \bar{C}_2 h^{\ell}, 
      \\ 
      \label{eq:tnormboundepb}
      &(\mu^b)^{-1/2}\tnorm{\boldsymbol{e}_{p^b}^{h,n}}_{q,b}
      \le \bar{C}_1 \Delta t + \bar{C}_2 h^{\ell},      
      \\
      \label{eq:tnormboundeps}
      &(\mu^s)^{-1/2}\big( \Delta t \sum_{i=1}^n \tnorm{\boldsymbol{e}_{p^s}^{h,i}}_{q,s}^2 \big)^{1/2}
      \le \bar{C}_1 \Delta t + \bar{C}_2 h^{\ell},      
    \end{align}
  \end{subequations}
  where the constants are defined as
  \small{ 
  \begin{equation*}
    \begin{split}
      \bar{C}_1
      =&
      C \bigg[ \max(c_0^{1/2},\lambda^{-1/2}\alpha) \norm[0]{\partial_{tt}p^p}_{L^1(J;L^2(\Omega^b))}
      + \lambda^{-1/2}\norm[0]{\partial_{tt}p^b}_{L^1(J;L^2(\Omega^b))}
      \\
      &\quad
      + \max(\gamma^{1/2}(\mu^s/\kappa)^{1/4}, \kappa^{-1/2})\norm[0]{\partial_{tt}u^{b}}_{L^2(J;L^2(\Gamma_I))}\bigg],
      \\
      \bar{C}_2
      =&
      C \bigg[\max\big((\mu^b)^{1/2}, (\mu^s  + \gamma(\mu^s/\kappa)^{1/2})^{1/2}\big) \norm[0]{\partial_tu^{b}}_{L^2(J;H^{\ell+1}(\Omega^b))}
      +\kappa^{-1/2} T^{1/2} \norm[0]{z}_{C^0(J;H^{\ell}(\Omega^b))}
      \\
      &\quad
      + \max\big((\mu^b)^{1/2}T, \del[1]{ \mu^s  + \gamma(\mu^s/\kappa)^{1/2}}^{1/2} T^{1/2} \big) \norm[0]{u^s}_{C^0(J;H^{\ell+1}(\Omega^s))}
      \bigg].
    \end{split}
  \end{equation*}
  }
\end{theorem}
\begin{proof}
  Choose $\boldsymbol{v}_h^s=\boldsymbol{e}_{u^s}^h$,
  $\boldsymbol{v}_h^b=d_t \boldsymbol{e}_{u^b}^h$,
  $\boldsymbol{q}_h^s=-\boldsymbol{e}_{p^s}^h$,
  $\boldsymbol{q}_h^b=-\boldsymbol{e}_{p^b}^h$,
  $\boldsymbol{q}^p_h=\boldsymbol{e}_{p^p}^h$, $w_h=e_{z}^h$ in
  \cref{eq:erroreqn} to find (we again suppress the superscript
  $n+1$):
  \begin{equation*}
    \begin{split}
      &a_h^s(\boldsymbol{e}_{u^s}^h, \boldsymbol{e}_{u^s}^h)  + a_h^b(\boldsymbol{e}_{u^b}^h, d_t\boldsymbol{e}_{u^b}^h) 
      + a_h^I((\bar{e}_{u^s}^h,d_t\bar{e}_{u^b}^h),(\bar{e}_{u^s}^h, d_t\bar{e}_{u^b}^h))
      \\
      &\qquad  + (\kappa^{-1} e_z^h, e_z^h)_{\Omega^b}
      + (c_0 d_te_{p^p}^h, e_{p^p}^h)_{\Omega^b}
      + c_h((d_te_{p^p}^h, d_te_{p^b}^h), \alpha e_{p^p}^h - e_{p^b}^h)
      \\
      &= a_h^I((0,\partial_tu^b-d_tu^b),(\bar{e}_{u^s}^h, d_t\bar{e}_{u^b}^h))
      + a_h^I((\bar{e}_{u^s}^I,d_t\bar{e}_{u^b}^I),(\bar{e}_{u^s}^h, d_t\bar{e}_{u^b}^h))
      + (\kappa^{-1} e_z^I, e_z^h)_{\Omega^b}
      \\
      &\qquad
      + (c_0 (\partial_t p^p - d_t p^p), e_{p^p}^h)_{\Omega^b}
      + c_h((\partial_t p^p - d_t p^p, \partial_t p^b-d_t p^b), \alpha e_{p^p}^h)
      \\
      &\qquad
      - b_h^I(\bar{e}_{p^p}^h,(0, \partial_t u^b - d_t u^b))
      - b_h^s(\boldsymbol{e}_{p^s}^h, \boldsymbol{e}_{u^s}^I)
      - b_h^b(\boldsymbol{e}_{p^b}^h, d_t\boldsymbol{e}_{u^b}^I).              
    \end{split}
  \end{equation*}
  Using the algebraic inequality $a(a-b) \ge (a^2-b^2)/2$ and
  multiplying the resulting equation by $\Delta t$, we obtain:
  \begin{equation}
    \label{eq:splitIs}
    \begin{split}
      \Delta t a_h^s&(\boldsymbol{e}_{u^s}^h, \boldsymbol{e}_{u^s}^h)
      + \gamma(\mu^s/\kappa)^{1/2}\Delta t\norm[0]{(\bar{e}_{u^s}^h - d_t\bar{e}_{u^b}^h)^t}_{\Gamma_I}^2
      + \kappa^{-1}\Delta t \norm[0]{e_z^h}_{\Omega^b}^2
      \\
      & + \tfrac{1}{2}\del[1]{a_h^b(\boldsymbol{e}_{u^b}^h, \boldsymbol{e}_{u^b}^h)
        - a_h^b(\boldsymbol{e}_{u^b}^{h,n}, \boldsymbol{e}_{u^b}^{h,n})}+ \tfrac{c_0}{2}\del[1]{\norm[0]{e_{p^p}^h}_{\Omega^b}^2 - \norm[0]{e_{p^p}^{h,n}}_{\Omega^b}^2}
      \\
      & + \tfrac{\lambda^{-1}}{2}\del[1]{\norm[0]{\alpha e_{p^p}^h - e_{p^b}^h}_{\Omega^b}^2
        - \norm[0]{\alpha e_{p^p}^{h,n} - e_{p^b}^{h,n}}_{\Omega^b}^2} 
      \\
      \le & \Delta t a_h^I((0,\partial_tu^b-d_tu^b),(\bar{e}_{u^s}^h, d_t\bar{e}_{u^b}^h))
      + \Delta t a_h^I((\bar{e}_{u^s}^I,d_t\bar{e}_{u^b}^I),(\bar{e}_{u^s}^h, d_t\bar{e}_{u^b}^h))
      \\
      &+ \Delta t (\kappa^{-1} e_z^I, e_z^h)_{\Omega^b}
      + c_0 \Delta t (\partial_t p^p - d_t p^p, e_{p^p}^h)_{\Omega^b}
      \\
      &+ \Delta t c_h((\partial_t p^p - d_t p^p, \partial_t p^b-d_t p^b), \alpha e_{p^p}^h)
      \\
      &
      - \Delta t b_h^I(\bar{e}_{p^p}^h,(0, \partial_t u^b - d_t u^b))
      - \Delta t b_h^s(\boldsymbol{e}_{p^s}^h, \boldsymbol{e}_{u^s}^I)
      - \Delta t b_h^b(\boldsymbol{e}_{p^b}^h, d_t\boldsymbol{e}_{u^b}^I)
      \\
      =:&I_1^{n+1} + \hdots + I_8^{n+1}.
    \end{split}
  \end{equation}
  Defining
  \begin{equation*}
    \begin{split}
      A_i^2
      &:= \tfrac{1}{2}a_h^b(\boldsymbol{e}_{u^b}^{h,i}, \boldsymbol{e}_{u^b}^{h,i})
      + \tfrac{c_0}{2}\norm[0]{e_{p^p}^{h,i}}_{\Omega^b}^2
      + \tfrac{\lambda^{-1}}{2}\norm[0]{\alpha e_{p^p}^{h,i} - e_{p^b}^{h,i}}_{\Omega^b}^2
      \\
      B_i^2
      &:= \tfrac{1}{2}\Delta t a_h^s(\boldsymbol{e}_{u^s}^{h,i}, \boldsymbol{e}_{u^s}^{h,i})
      + \tfrac{\kappa^{-1}}{2}\Delta t \norm[0]{e_z^{h,i}}_{\Omega^b}^2
      + \tfrac{\gamma(\mu^s/\kappa)^{1/2}}{2}\Delta t\norm[0]{(\bar{e}_{u^s}^{h,i} - d_t\bar{e}_{u^b}^{h,i})^t}_{\Gamma_I}^2,
    \end{split}
  \end{equation*}
  and by \cref{eq:initialconditions0}, we can write \cref{eq:splitIs}
  as
  \begin{equation}
    \label{eq:splitIs-compact}
    A_{n+1}^2 + 2 B_{n+1}^2 \le A_n^2 + \sum_{k=1}^8I_k^{n+1}, \qquad A_0 = 0.
  \end{equation}
  We will now bound each of the terms $I_k^{n+1}$, $k=1,\hdots,8$.

  First, by the Cauchy--Schwarz and Young's inequalities, we find that for any $\delta > 0$
  \begin{equation*}
    \begin{split}
      I_1^{n+1}
      &\le \Delta t \gamma(\mu^s/\kappa)^{1/2}\norm[0]{(\partial_tu^b-d_tu^b)^t}_{\Gamma_I}\norm[0]{(\bar{e}_{u^s}^h-d_t\bar{e}_{u^b}^h)^t}_{\Gamma_I}
      \\
      &\le \Delta t \frac{\gamma(\mu^s/\kappa)^{1/2}}{2\delta}\norm[0]{(\partial_tu^b-d_tu^b)^t}_{\Gamma_I}^2
      + \delta \Delta t \frac{\gamma(\mu^s/\kappa)^{1/2}}{2} \norm[0]{(\bar{e}_{u^s}^h-d_t\bar{e}_{u^b}^h)^t}_{\Gamma_I}^2
      \\
      &\le \Delta t \frac{\gamma(\mu^s/\kappa)^{1/2}}{2\delta}\norm[0]{(\partial_tu^b-d_tu^b)^t}_{\Gamma_I}^2
      + \delta B_{n+1}^2.
    \end{split}
  \end{equation*}
  Similarly,
  \begin{equation*}
    \begin{split}
      I_2^{n+1}
      &\le \Delta t \frac{\gamma(\mu^s/\kappa)^{1/2}}{2\delta}\norm[0]{(\bar{e}_{u^s}^I - d_t\bar{e}_{u^b}^I)^t}_{\Gamma_I}^2
      + \delta B_{n+1}^2,
      \\
      I_3^{n+1}
      &\le \Delta t \frac{\kappa^{-1}}{2\delta}\norm[0]{e_z^I}_{\Omega^b}^2 + \delta B_{n+1}^2,
    \end{split}
  \end{equation*}
  while application of the Cauchy--Schwarz inequality results in
  \begin{equation*}
    I_4^{n+1}
    \le c_0\Delta t \norm[0]{\partial_tp^p - d_tp^p}_{\Omega^b}\norm[0]{e_{p^p}^h}_{\Omega^b}
    \le (2c_0)^{1/2}\Delta t \norm[0]{\partial_tp^p - d_tp^p}_{\Omega^b} A_{n+1}.
  \end{equation*}
  For $I_5^{n+1}$ we find, using the Cauchy--Schwarz and triangle
  inequalities and $C_* \mu^b \le \lambda$, that
  \begin{equation*}
    \begin{split}
      I_5^{n+1}
      &\le \Delta t \lambda^{-1}\norm[0]{\alpha(\partial_tp^p-d_tp^p) - (\partial_tp^b-d_tp^b)}_{\Omega^b}\norm[0]{\alpha e_{p^p}^h}_{\Omega^b}
      \\
      &\le \Delta t \lambda^{-1}\del[1]{\norm[0]{\alpha(\partial_tp^p-d_tp^p)}_{\Omega^b} + \norm[0]{\partial_tp^b-d_tp^b}_{\Omega^b}}
      \del[1]{\norm[0]{\alpha e_{p^p}^h - e_{p^b}^h}_{\Omega^b} + \norm[0]{e_{p^b}^h}_{\Omega^b}}
      \\
      &\le \Delta t \lambda^{-1/2}\del[1]{\norm[0]{\alpha(\partial_tp^p-d_tp^p)}_{\Omega^b} + \norm[0]{\partial_tp^b-d_tp^b}_{\Omega^b}} \\
      &\quad \times \del[1]{\lambda^{-1/2}\norm[0]{\alpha e_{p^p}^h - e_{p^b}^h}_{\Omega^b} + (C_*\mu^b)^{-1/2}\norm[0]{e_{p^b}^h}_{\Omega^b}}.      
    \end{split}
  \end{equation*}
  To bound this further, let $\boldsymbol{v}_h^s = \boldsymbol{0}$ and
  $\boldsymbol{v}_h^b \in \widetilde{V}_h^b$ in \cref{eq:erroreqn-a}
  to find that
  $b_h^b(\boldsymbol{e}_{p^b}^h, \boldsymbol{v}_h^b) =
  -a_h^b(\boldsymbol{e}_{u^b}^h, \boldsymbol{v}_h^b)$. By
  \cref{eq:ah-cont-j,eq:ah-coer,eq:total-pressure-inf-sup-alt},
  \begin{equation}
    \label{eq:boundqhqbinfsups}
    \tnorm{\boldsymbol{e}_{p^b}^{h}}_{q,b} \le
    C\sup_{\boldsymbol{0} \ne \boldsymbol{v}_h \in \widetilde{\boldsymbol{V}}_h^b}
    \frac{b_h^b(\boldsymbol{e}_{p^b}^{h}, \boldsymbol{v}_h)}{\tnorm{\boldsymbol{v}_h}_{v,b}}
    \le C\mu^b\tnorm{\boldsymbol{e}_{u^b}^h}_{v,b} \le C (\mu^b)^{1/2}a_h^b(\boldsymbol{e}_{u^b}^h,\boldsymbol{e}_{u^b}^h)^{1/2}.
  \end{equation}
  In other words,
  $(\mu^b)^{-1/2}\norm[0]{e_{p^b}^h}_{\Omega^b} \le
  Ca_h^b(\boldsymbol{e}_{u^b}^h,\boldsymbol{e}_{u^b}^h)^{1/2} \le
  CA_{n+1}$ so that
  \begin{equation*}
    I_5^{n+1}
    \le C\Delta t \lambda^{-1/2}\del[1]{\alpha\norm[0]{\partial_tp^p-d_tp^p}_{\Omega^b} + \norm[0]{\partial_tp^b-d_tp^b}_{\Omega^b}}
    A_{n+1}.
  \end{equation*}
  To bound $I_6^{n+1}$ we use \cref{lem:pPGammaI_z} and
  Cauchy--Schwarz, triangle, and Young's inequalities:
  \begin{equation*}
    \begin{split}
      I_6^{n+1}
      =&
      \Delta t \langle \bar{e}_{p^p}^h, (\partial_tu^b-d_tu^b) \cdot n \rangle_{\Gamma_I}
      \le C\kappa^{-1}\Delta t\norm[0]{z - z_h}_{\Omega^b} \norm[0]{d_tu^b - \partial_tu^b}_{\Gamma_I}
      \\
      \le& C\kappa^{-1}\Delta t \norm[0]{e_z^h}_{\Omega^b} \norm[0]{d_tu^b - \partial_tu^b}_{\Gamma_I}
      + C\kappa^{-1}\Delta t \norm[0]{e_z^I}_{\Omega^b} \norm[0]{d_tu^b - \partial_tu^b}_{\Gamma_I}
      \\
      \le& \delta B_{n+1}^2
      + \frac{C\kappa^{-1}}{2\delta} \Delta t \norm[0]{d_tu^b - \partial_tu^b}_{\Gamma_I}^2
      + \frac{C\kappa^{-1}}{2}\Delta t \norm[0]{e_z^I}_{\Omega^b}^2
      + \frac{C\kappa^{-1}}{2}\Delta t \norm[0]{d_tu^b - \partial_tu^b}_{\Gamma_I}^2.
    \end{split}
  \end{equation*}
  To bound $I_7^{n+1}$ and $I_8^{n+1}$ we first bound $\tnorm{\boldsymbol{e}_p^h}_q$. Observe from
  \cref{eq:erroreqn-a,eq:ah-cont-j}, and the Cauchy--Schwarz inequality
  that for all $\boldsymbol{v}_h \in \widehat{\boldsymbol{V}}_h$, since $b_h^I(\bar{e}_{p^p}^h,(\bar{v}_h^s,\bar{v}_h^b))=0$,
  \begin{equation*}
    \begin{split}
      b_h(\boldsymbol{e}_p^h, \boldsymbol{v}_h)
      =&
      - a_h(\boldsymbol{e}_u^h, \boldsymbol{v}_h)
      - a_h^I((\bar{e}_{u^s}^h,d_t\bar{e}_{u^b}^h),(\bar{v}_h^s,\bar{v}_h^b))
      \\
      &+ a_h^I((0,\partial_tu^b-d_tu^b),(\bar{v}_h^s,\bar{v}_h^b))
      + a_h^I((\bar{e}_{u^s}^I,d_t\bar{e}_{u^b}^I),(\bar{v}_h^s,\bar{v}_h^b))
      \\
      \le& C\del[1]{ \mu^s \tnorm{\boldsymbol{e}_{u^s}^h}_{v,s}\tnorm{\boldsymbol{v}_h^s}_{v,s}
        + \mu^b \tnorm{\boldsymbol{e}_{u^b}^h}_{v,b}\tnorm{\boldsymbol{v}_h^b}_{v,b} }
      \\
      &+ \gamma(u^s/\kappa)^{1/2}\norm[0]{(\bar{e}_{u^s}^h - d_t\bar{e}_{u^b}^h)^t}_{\Gamma_I}\tnorm{\boldsymbol{v}_h}_v
      \\
      &+ \gamma(u^s/\kappa)^{1/2}\del[1]{ \norm[0]{(\partial_tu^b-d_tu^b)^t}_{\Gamma_I}
        + \norm[0]{(\bar{e}_{u^s}^I - d_t\bar{e}_{u^b}^I)^t}_{\Gamma_I} }\tnorm{\boldsymbol{v}_h}_v.
    \end{split}
  \end{equation*}
  By \cref{eq:inf-sup,eq:ah-coer} we then find that
  \begin{equation}
    \begin{split}
      \tnorm{\boldsymbol{e}_p^h}_q
      \le &
      C \big\{
      \mu^s \tnorm{\boldsymbol{e}_{u^s}^h}_{v,s}
      + \mu^b \tnorm{\boldsymbol{e}_{u^b}^h}_{v,b}
      + \gamma(u^s/\kappa)^{1/2}\norm[0]{(\bar{e}_{u^s}^h - d_t\bar{e}_{u^b}^h)^t}_{\Gamma_I}
      \\
      &+ \gamma(u^s/\kappa)^{1/2}\del[1]{ \norm[0]{(\partial_tu^b-d_tu^b)^t}_{\Gamma_I}
        + \norm[0]{(\bar{e}_{u^s}^I - d_t\bar{e}_{u^b}^I)^t}_{\Gamma_I} } \big\}
      \\
      \le &
      C \big\{ (\Delta t)^{-1/2} \sbr[1]{ (\mu^s)^{1/2}  + \gamma^{1/2}(\mu^s/\kappa)^{1/4}} B_{n+1}
      + (\mu^b)^{1/2} A_{n+1}
      \\
      &+ \gamma(u^s/\kappa)^{1/2}\del[1]{ \norm[0]{(\partial_tu^b-d_tu^b)^t}_{\Gamma_I}
        + \norm[0]{(\bar{e}_{u^s}^I - d_t\bar{e}_{u^b}^I)^t}_{\Gamma_I} } \big\}.      
    \end{split}
  \end{equation}
  Therefore, using \cref{ineq:bhj-bound},
  \begin{equation*}
    \begin{split}
      I_7^{n+1}
      \le
      & C \Delta t \tnorm{\boldsymbol{e}_{p^s}^h}_{q,s} \tnorm{\boldsymbol{e}_{u^s}^I}_{v,s}
      \\
      \le
      & C\Delta t \cbr[1]{ (\Delta t)^{-1/2} \sbr[1]{ (\mu^s)^{1/2}  + \gamma^{1/2}(\mu^s/\kappa)^{1/4}} B_{n+1}
        + (\mu^b)^{1/2} A_{n+1}} \tnorm{\boldsymbol{e}_{u^s}^I}_{v,s}
      \\
      & + C\Delta t \cbr[1]{ \gamma(u^s/\kappa)^{1/2}\del[1]{ \norm[0]{(\partial_tu^b-d_tu^b)^t}_{\Gamma_I}
          + \norm[0]{(\bar{e}_{u^s}^I - d_t\bar{e}_{u^b}^I)^t}_{\Gamma_I} } } \tnorm{\boldsymbol{e}_{u^s}^I}_{v,s}
      \\
      =: & I_{71} + I_{72}.
    \end{split}    
  \end{equation*}
  Young's inequality is used to bound $I_{71}$:
  \begin{equation*}
    \begin{split}
      I_{71}
      =&
      C (\Delta t)^{1/2} \sbr[1]{ (\mu^s)^{1/2}  + \gamma^{1/2}(\mu^s/\kappa)^{1/4}} B_{n+1}\tnorm{\boldsymbol{e}_{u^s}^I}_{v,s}
      + C\Delta t(\mu^b)^{1/2} A_{n+1}\tnorm{\boldsymbol{e}_{u^s}^I}_{v,s}
      \\
      \le& \delta B_{n+1}^2
      + C \Delta t \del[1]{ \mu^s  + \gamma(\mu^s/\kappa)^{1/2}}\tnorm{\boldsymbol{e}_{u^s}^I}_{v,s}^2
      + C\Delta t(\mu^b)^{1/2} A_{n+1}\tnorm{\boldsymbol{e}_{u^s}^I}_{v,s}.
    \end{split}
  \end{equation*}
  Similarly, we find that
  \begin{equation*}
    \begin{split}
      I_8^{n+1}
      \le
      & C \Delta t \tnorm{\boldsymbol{e}_{p^b}^h}_{q,b} \tnorm{d_t\boldsymbol{e}_{u^b}^I}_{v,b}
      \\
      \le
      & C\Delta t \cbr[1]{ (\Delta t)^{-1/2} \sbr[1]{ (\mu^s)^{1/2}  + \gamma^{1/2}(\mu^s/\kappa)^{1/4}} B_{n+1}
        + (\mu^b)^{1/2} A_{n+1}} \tnorm{d_t\boldsymbol{e}_{u^b}^I}_{v,b}
      \\
      & + C\Delta t \cbr[1]{ \gamma(u^s/\kappa)^{1/2}\del[1]{ \norm[0]{(\partial_tu^b-d_tu^b)^t}_{\Gamma_I}
          + \norm[0]{(\bar{e}_{u^s}^I - d_t\bar{e}_{u^b}^I)^t}_{\Gamma_I} } } \tnorm{d_t\boldsymbol{e}_{u^b}^I}_{v,b}
      \\
      =: & I_{81} + I_{82},
    \end{split}    
  \end{equation*}
  and
  \begin{equation*}
    I_{81} \le \delta B_{n+1}^2
      + C \Delta t \del[1]{ \mu^s  + \gamma(\mu^s/\kappa)^{1/2}}\tnorm{d_t\boldsymbol{e}_{u^b}^I}_{v,b}^2
      + C\Delta t(\mu^b)^{1/2} A_{n+1}\tnorm{d_t\boldsymbol{e}_{u^b}^I}_{v,b}.
  \end{equation*}
  Adding up the various bounds for $I_k^{n+1}$, $k=1,\hdots,8$, we
  find
  \begin{equation}
    \label{eq:boundsSumIkall}
    \sum_{k=1}^{8} I_k^{n+1} \le 6\delta B_{n+1}^2 + E_{n+1}A_{n+1} + D_{n+1},
  \end{equation}
  where
  \begin{equation*}
    \begin{split}
      E_i
      = C\Delta t \Big[& c_0^{1/2}\norm[0]{\partial_tp^{p,i} - d_tp^{p,i}}_{\Omega^b}
      + \lambda^{-1/2}\alpha\norm[0]{\partial_tp^{p,i}-d_tp^{p,i}}_{\Omega^b}
      \\
      &+ \lambda^{-1/2}\norm[0]{\partial_tp^{b,i}-d_tp^{b,i}}_{\Omega^b}
      + (\mu^b)^{1/2} \tnorm{\boldsymbol{e}_{u^s}^{I,i}}_{v,s}
      + (\mu^b)^{1/2}\tnorm{d_t\boldsymbol{e}_{u^b}^{I,i}}_{v,b}\Big],
    \end{split}
  \end{equation*}
  and
  \begin{equation*}
    \begin{split}
      D_i
      =
     C\Delta t \Big[ & \gamma(\mu^s/\kappa)^{1/2}\norm[0]{(\partial_tu^{b,i}-d_tu^{b,i})^t}_{\Gamma_I}^2
      + \gamma(\mu^s/\kappa)^{1/2}{2\delta}\norm[0]{(\bar{e}_{u^s}^{I,i} - d_t\bar{e}_{u^b}^{I,i})^t}_{\Gamma_I}^2
      \\
      &+ \kappa^{-1}\norm[0]{e_z^{I,i}}_{\Omega^b}^2
      + \kappa^{-1} \norm[0]{d_tu^{b,i} - \partial_tu^{b,i}}_{\Gamma_I}^2 
      + \del[1]{ \mu^s  + \gamma(\mu^s/\kappa)^{1/2}}\tnorm{\boldsymbol{e}_{u^s}^{I,i}}_{v,s}^2
      \\
      &+ \gamma(\mu^s/\kappa)^{1/2}\del[1]{ \norm[0]{(\partial_tu^{b,i}-d_tu^{b,i})^t}_{\Gamma_I}
          + \norm[0]{(\bar{e}_{u^s}^{I,i} - d_t\bar{e}_{u^b}^{I,i})^t}_{\Gamma_I} } \tnorm{\boldsymbol{e}_{u^s}^{I,i}}_{v,s}
      \\
      &+ \gamma(\mu^s/\kappa)^{1/2}\del[1]{ \norm[0]{(\partial_tu^{b,i}-d_tu^{b,i})^t}_{\Gamma_I}
          + \norm[0]{(\bar{e}_{u^s}^{I,i} - d_t\bar{e}_{u^b}^{I,i})^t}_{\Gamma_I} } \tnorm{d_t\boldsymbol{e}_{u^b}^{I,i}}_{v,b}
      \\
      &+ \del[1]{ \mu^s  + \gamma(\mu^s/\kappa)^{1/2}}\tnorm{d_t\boldsymbol{e}_{u^b}^{I,i}}_{v,b}^2\Big].
    \end{split}
  \end{equation*}
  Combining \cref{eq:splitIs-compact,eq:boundsSumIkall} (with
  $\delta = 1/6$) and summing over the time levels, we obtain
  \begin{equation*}
    A_n^2 + \sum_{i=1}^nB_i^2 \le \sum_{i=1}^n E_iA_i + \sum_{i=1}^n D_i.
  \end{equation*}
  By \cite[Lemma 4.1]{Cesmelioglu:2022},
  \begin{equation}
    \label{eq:splitIs-compact-full}
    A_n + \big( \sum_{i=1}^n B_i^2 \big)^{1/2} \le C \del[2]{ \sum_{i=1}^nE_i + \big( \sum_{i=1}^n D_i \big)^{1/2} }.
  \end{equation}
  We will now bound the two sums on the right hand side
  separately. For this we require the following inequalities (that can
  be proven by Taylor series expansions) \cite[Lemma 3.2]{Cesmelioglu:2020a}, \cite[(4.5a)]{Cesmelioglu:2022}:
  \begin{subequations}
    \begin{align}
      \label{ineq:taylor1}
      \sum_{i=1}^n \Delta t  \|\partial_t u^{b,i}-d_t u^{b,i}\|^2_{\Gamma_I}
      &\le C(\Delta t)^2 \| \partial_{tt} u^b \|_{L^2(J; L^2(\Gamma_I))}, && 
      \\
      \label{ineq:taylor2}
      \sum_{i=1}^n \Delta t  \|\partial_t \psi^i-d_t\psi^i\|_{\Omega^b}
      &\le \Delta t \| \partial_{tt} \psi \|_{L^1(J; L^2(\Omega^b))}, && \psi = p^b, p^p,
      \\
      \label{ineq:taylor3}
      \Delta t \sum_{i=1}^n \| d_t u^{b,i} \|_{\ell+1,\Omega^j} 
      &\le C \| \partial_t u^b \|_{L^1(J; H^{\ell+1}(\Omega^b))}, && 
      \\
      \label{ineq:taylor3b}
      \Delta t \sum_{i=1}^n \| d_t u^{b,i} \|^2_{\ell+1,\Omega^b} 
      &\le C \| \partial_t u^b \|^2_{L^2(J; H^{\ell+1}(\Omega^b))}, &&
      \\
      \label{ineq:taylor4}
      \Delta t \sum_{i=1}^n \| d_t \bar{e}_{u^b}^{I,i}\|^2_{\Gamma_I} 
      &\le C h^{2\ell} \| \partial_t u^b \|^2_{L^2(J; H^{\ell+1}(\Omega^b))}, &&
    \end{align}
  \end{subequations}
  and the inequalities (that can be proven using
  \cref{eq:ah-coer,eq:interpolation-property-ell,ineq:taylor3,ineq:taylor3b})
  \begin{equation*}
    \begin{split}
      \Delta t \sum_{i=1}^n \tnorm{d_t\boldsymbol{e}_{u^b}^{I,i}}_{v,b}
      &\leq C \Delta t (\mu^b)^{-1/2}\sum_{i=1}^n a_h^b(d_t\boldsymbol{e}_{u^b}^{I,i}, d_t\boldsymbol{e}_{u^b}^{I,i})^{1/2}
      \\
      &\leq C \Delta t \sum_{i=1}^n h^{\ell}\norm[0]{d_t u^{b,i}}_{\ell+1,\Omega^b}
      \leq Ch^{\ell} \norm[0]{\partial_t u^{b}}_{L^1(J;H^{\ell+1}(\Omega^b))},
      \\
      \Delta t \sum_{i=1}^n \tnorm{d_t\boldsymbol{e}_{u^b}^{I,i}}_{v,b}^2
      &\leq C \Delta t (\mu^{b})^{-1}\sum_{i=1}^n a_h^b(d_t\boldsymbol{e}_{u^b}^{I,i}, d_t\boldsymbol{e}_{u^b}^{I,i})
      \\
      &\leq C \Delta t \sum_{i=1}^n h^{2\ell} \norm[0]{d_tu^{b,i}}^2_{\ell+1,\Omega^b}
      \leq Ch^{2\ell} \norm[0]{\partial_tu^{b}}^2_{L^2(J;H^{\ell+1}(\Omega^b))},
      \\
      \Delta t \sum_{i=1}^n \tnorm{\boldsymbol{e}_{u^s}^{I,i}}_{v,s}
      &\leq C \Delta t (\mu^s)^{-1/2} \sum_{i=1}^n a_h^s(\boldsymbol{e}_{u^s}^{I,i},\boldsymbol{e}_{u^s}^{I,i})^{1/2}
      \\
      &\leq C \Delta t \sum_{i=1}^n h^{\ell} \norm[0]{u^{s,i}}_{\ell+1,\Omega^s}
      \le C T h^{\ell} \norm[0]{u^{s}}_{C^0(J;H^{\ell+1}(\Omega^s))}.
    \end{split}
  \end{equation*}
  Similarly, we have by \cref{eq:interpolation-property-bdm,eq:interpolation-property-ell},
  \begin{equation*}
    \begin{split}
      \Delta t \sum_{i=1}^n \norm[0]{e_z^{I,i}}_{\Omega^b}^2 &\le C T h^{2\ell} \norm[0]{z}_{C^0(J;H^{\ell}(\Omega^b))}^2,
      \\
      \Delta t \sum_{i=1}^n \norm[0]{\bar{e}_{u^s}^{I,i}}_{\Gamma_I}^2 &\le C T h^{2\ell} \norm[0]{u^s}_{C^0(J;H^{\ell+1}(\Omega^s))}^2. 
    \end{split}
  \end{equation*}
  We now find:
  \begin{equation*}
    \begin{split}
      \sum_{i=1}^nE_i
      \le&
      C \Delta t\bigg[ c_0^{1/2} \norm[0]{\partial_{tt}p^p}_{L^1(J;L^2(\Omega^b))}
      + \lambda^{-1/2}\alpha\norm[0]{\partial_{tt}p^p}_{L^1(J;L^2(\Omega^b))}
       \\
       &\qquad
      + \lambda^{-1/2}\norm[0]{\partial_{tt}p^b}_{L^1(J;L^2(\Omega^b))}\bigg]
      \\
      &+ C h^{\ell}\bigg[ (\mu^b)^{1/2}T \norm[0]{u^{s}}_{C^0(J;H^{\ell+1}(\Omega^s))}
      + (\mu^b)^{1/2} \norm[0]{\partial_tu^{b}}_{L^2(J;H^{\ell+1}(\Omega^b))}\bigg],
    \end{split}
  \end{equation*}
  and, after applying Young's inequality to $D_i$, 
  \begin{equation*}
    \begin{split}
      \bigg(\sum_{i=1}^nD_i\bigg)^{1/2}
      \le
      & C \Delta t \bigg[ \gamma^{1/2}(\mu^s/\kappa)^{1/4}\norm[0]{\partial_{tt}u^{b}}_{L^2(J;L^2(\Gamma_I))}
       + \kappa^{-1/2} \norm[0]{\partial_{tt}u^{b}}_{L^2(J;L^2(\Gamma_I))} \bigg]
      \\
      &+ C h^{\ell} \bigg[\kappa^{-1/2} T^{1/2} \norm[0]{z}_{C^0(J;H^{\ell}(\Omega^b))}
       \\
       &\qquad\qquad
      + \del[1]{ \mu^s  + \gamma(\mu^s/\kappa)^{1/2}}^{1/2} T^{1/2} \norm[0]{u^s}_{C^0(J;H^{\ell+1}(\Omega^s))}
      \\
      &\qquad\qquad
      + \del[1]{ \mu^s  + \gamma(\mu^s/\kappa)^{1/2}}^{1/2} \norm[0]{\partial_tu^b}_{L^2(J;H^{\ell+1}(\Omega^b))} \bigg],
    \end{split}
  \end{equation*}
  From \cref{eq:splitIs-compact-full} we then find:
  \small{ 
  \begin{equation*}
    \begin{split}
      &A_n + \big( \sum_{i=1}^n B_i^2 \big)^{1/2}
      \\
      &\le
      C \Delta t\bigg[ \max(c_0^{1/2},\lambda^{-1/2}\alpha) \norm[0]{\partial_{tt}p^p}_{L^1(J;L^2(\Omega^b))}
      + \lambda^{-1/2}\norm[0]{\partial_{tt}p^b}_{L^1(J;L^2(\Omega^b))}
      \\
      &\qquad\qquad
      + \max(\gamma^{1/2}(\mu^s/\kappa)^{1/4}, \kappa^{-1/2})\norm[0]{\partial_{tt}u^{b}}_{L^2(J;L^2(\Gamma_I))}\bigg]
      \\
      &+ C h^{\ell}\bigg[\max\big((\mu^b)^{1/2}, (\mu^s  + \gamma(\mu^s/\kappa)^{1/2})^{1/2}\big) \norm[0]{\partial_tu^{b}}_{L^2(J;H^{\ell+1}(\Omega^b))}
      +\kappa^{-1/2} T^{1/2} \norm[0]{z}_{C^0(J;H^{\ell}(\Omega^b))}
      \\
      &\qquad\qquad+ \max\big((\mu^b)^{1/2}T, \del[1]{ \mu^s  + \gamma(\mu^s/\kappa)^{1/2}}^{1/2} T^{1/2} \big) \norm[0]{u^s}_{C^0(J;H^{\ell+1}(\Omega^s))}
      \bigg].     
    \end{split}
  \end{equation*}
  }
  \Cref{ineq:apriori-1} now follows by definition of $A_i$ and $B_i$
  and the coercivity of $a_h^s$ and $a_h^b$ \cref{eq:ah-coer} while
  \cref{eq:tnormboundepb} follows from \cref{eq:boundqhqbinfsups} and
  noting that
  $a_h^b(\boldsymbol{e}_{u^b}^{h,n}, \boldsymbol{e}_{u^b}^{h,n})^{1/2}
  \le A_n$. Finally, \cref{eq:tnormboundeps} follows using similar
  steps as used to find \cref{eq:tnormboundepb}: let
  $\boldsymbol{v}_h^b = \boldsymbol{0}$ and
  $\boldsymbol{v}_h^s \in \widetilde{V}_h^s$ in \cref{eq:erroreqn-a}
  to find that
  $b_h^s(\boldsymbol{e}_{p^s}^{h,i}, \boldsymbol{v}_h^s) =
  -a_h^s(\boldsymbol{e}_{u^s}^{h,i}, \boldsymbol{v}_h^s)$. By
  \cref{eq:ah-cont-j,eq:ah-coer,eq:total-pressure-inf-sup-alt},
  \begin{equation*}
    C_1 \tnorm{\boldsymbol{e}_{p^s}^{h,i}}_{q,s} \le
    \sup_{\boldsymbol{0} \ne \boldsymbol{v}_h \in \widetilde{\boldsymbol{V}}_h^s}
    \frac{b_h^s(\boldsymbol{e}_{p^s}^{h,i}, \boldsymbol{v}_h)}{\tnorm{\boldsymbol{v}_h}_{v,s}}
    \le C_2\mu^s\tnorm{\boldsymbol{e}_{u^s}^{h,i}}_{v,s}
    \le C_3 (\mu^s)^{1/2}a_h^s(\boldsymbol{e}_{u^s}^{h,i},\boldsymbol{e}_{u^s}^{h,i})^{1/2}.
  \end{equation*}
  The result follows noting that
  $(\mu^s)^{-1} \Delta t \tnorm{\boldsymbol{e}_{p^s}^{h,i}}_{q,s}^2
  \le C B_i^2$. 
\end{proof}

The main result of this section, an a priori error estimate for the
solution to the HDG method that is robust in the limits
$\lambda\to\infty$ and $c_0\to 0$, is now a consequence of
\cref{thm:apriorierrorestimate}.
\begin{corollary}
  Let $(u,p,z,p^p)$ be the solution to the coupled Stokes--Biot problem
  \cref{eq:stokesbiot,eq:interface,eq:bics} on time interval $J=(0,T]$.
  In addition to the regularity assumptions used in
  \cref{thm:apriorierrorestimate} we further assume that $p^s \in C^0(J;H^{\ell}(\Omega^s))$. Let
  $(\boldsymbol{u}_h^{n}, \boldsymbol{p}_h^{n}, z_h^{n},
  \boldsymbol{p}_h^{p,n}) \in \boldsymbol{X}_h$, for $n=1,\hdots,N$,
  be the solution to \cref{eq:fullydiscreteHDG} with initial
  conditions \cref{eq:initialconditions0}. Define
  $\boldsymbol{u}^{j} = (u^{j},\gamma_{V^j}(u^j))$. Then
  \begin{subequations}
  \begin{equation}
      \begin{split}
      \label{ineq:apriori-1-f}
      &(\mu^b)^{1/2}\tnorm{\boldsymbol{u}^{b,n} - \boldsymbol{u}_h^{b,n}}_{v,b}
        + c_0^{1/2} \norm[0]{p^{p,n}-p_h^{p,n}}_{\Omega^b}
      +\lambda^{-1/2}\norm[0]{\alpha(p^{p,n}-p_h^{p,n}) - (p^{b,n} - p^{b,n}_h)}_{\Omega^b}
      \\
      & \qquad +(\mu^s)^{1/2}\big(\Delta t \sum_{i=1}^n \tnorm{\boldsymbol{u}^{s,i} - \boldsymbol{u}^{s,i}_h}_{v,s}^2\big)^{1/2} 
        + \kappa^{-1/2} \big(\Delta t\sum_{i=1}^n\norm[0]{z^i - z_h^i}^2_{\Omega^b}\big)^{1/2}
      \\
      & \qquad + \gamma^{1/2}(\mu^s/\kappa)^{1/4} \big( \Delta t \sum_{i=1}^n \norm[0]{( (u_s^i - \bar{u}_h^{s,i}) - d_t(u_b^i - \bar{u}_h^{b,i}))^t}_{\Gamma_I}^2\big)^{1/2}
      \leq \bar{C}_1\Delta t + \tilde{C}_2 h^{\ell},          
      \end{split}
  \end{equation}
  \begin{equation}
      \label{eq:tnormboundepb-f}
      (\mu^b)^{-1/2}\tnorm{\boldsymbol{p}^n - \boldsymbol{p}_h^{b,n}}_{q,b}
      \le \bar{C}_1 \Delta t + (\bar{C}_2 + C (\mu^b)^{-1/2}\norm[0]{ p^{b} }_{C^0(J;H^{\ell}(\Omega^b))}) h^{\ell},      
  \end{equation}
  \begin{equation}
      \label{eq:tnormboundeps-f}
      (\mu^s)^{-1/2}\big(\Delta t \sum_{i=1}^n \tnorm{\boldsymbol{p}^{s,i}- \boldsymbol{p}_h^{s,i}}_{q,s}^2\big)^{1/2}
      \le 2\bar{C}_1\Delta t + \widehat{C}_2h^{\ell},
  \end{equation}
  \end{subequations}
  where $\bar{C}_1$ and $\bar{C}_2$ are the constants defined in
  \cref{thm:apriorierrorestimate} and
  \begin{equation*}
    \begin{split}
      \tilde{C}_2
      =&
      \bar{C}_2 + C \bigg[ (\mu^b)^{1/2}\|u^b\|_{C^0(J;H^{\ell+1}(\Omega^b))}
      + \max((\mu^s)^{1/2}, \gamma^{1/2}(\mu^s/\kappa)^{1/4})T^{1/2}\|u^s\|_{C^0(J;H^{\ell+1}(\Omega^s))}
      \\
      &\qquad\quad
      + \gamma^{1/2}(\mu^s/\kappa)^{1/4}\|\partial_tu^b\|_{L^2(J;H^{\ell+1}(\Omega^b)}
      + \max(c_0^{1/2}, \lambda^{-1/2}\alpha)\|p^p\|_{C^0(J;H^{\ell}(\Omega^b))}
      \\
      &\qquad\quad      
      + \lambda^{-1/2}\|p^b\|_{C^0(J;H^{\ell}(\Omega^b))}
      + \kappa^{-1/2}T^{1/2}\|z\|_{C^0(J;H^{\ell}(\Omega^b))}
      \bigg],
      \\
      \widehat{C}_2
      =& \del[1]{4\bar{C}_2^2 + C(\mu^s)^{-1}T \norm[0]{ p^{s} }_{C^0(J;H^{\ell}(\Omega^s))}^2}^{1/2}.
    \end{split}
  \end{equation*}
\end{corollary}
\begin{proof}
  \Cref{ineq:apriori-1-f} is a direct consequence of the triangle
  inequality, \cref{ineq:apriori-1}, and the following estimates:
  \begin{align*}
    \tnorm{\boldsymbol{e}_{u^b}^{I,n}}_{v,b} & \leq Ch^{\ell}\|u^b\|_{C^0(J;H^{\ell+1}(\Omega^b))},
    \\
    \|e_{p^p}^{I,n}\|_{\Omega^b}&\leq Ch^{\ell}\|p^p\|_{C^0(J;H^{\ell}(\Omega^b))},
    \\
    \|\alpha e_{p^p}^{I,n}-e_{p^b}^{I,n}\|_{\Omega^b}&\leq C h^{\ell}(\alpha\|p^p\|_{C^0(J;H^{\ell}(\Omega^b))}+\|p^b\|_{C^0(J;H^{\ell}(\Omega^b))}),
    \\
    \Big(\Delta t \sum_{i=1}^n\tnorm{\boldsymbol{e}_{u^s}^{I,i}}_{v,s}^2\Big)^{1/2} & \leq CT^{1/2} h^{\ell}\|u^s\|_{C^0(J;H^{\ell+1}(\Omega^s))},
    \\
    \Big(\Delta t \sum_{i=1}^n \|e_{z}^{I,i}\|_{\Omega^b}^2\Big)^{1/2}&\leq CT^{1/2} h^{\ell}\|z\|_{C^0(J;H^{\ell}(\Omega^b))},
    \\
    \Big(\Delta t \sum_{i=1}^n \|\bar{e}_{u^s}^{I,i}-d_t\bar{e}_{u^b}^{I,i}\|_{\Gamma_I}^2\Big)^{1/2}&\leq Ch^{\ell}\Big( T^{1/2} \|u^s\|_{C^0(J;H^{\ell+1}(\Omega^s))}+\|\partial_tu^b\|_{L^2(J;H^{\ell+1}(\Omega^b))}\Big).
  \end{align*}
  Next, note that
  \begin{equation}
    \label{eq:epjinterpolationbound}
    \tnorm{ \boldsymbol{e}_{p^j}^{I,i} }_{q,j} \le C h^{\ell} \norm[0]{ p^{j,i} }_{\ell, \Omega^j}, \qquad j=s,b.
  \end{equation}
  \Cref{eq:tnormboundepb-f} follows by a triangle inequality,
  \cref{eq:tnormboundepb}, and
  \cref{eq:epjinterpolationbound}. Finally, to show
  \cref{eq:tnormboundeps-f} we note that, by a triangle inequality,
  \cref{eq:epjinterpolationbound}, and Young's inequality
  \begin{equation*}
    \begin{split}
      (\mu^s)^{-1} \Delta t \sum_{i=1}^n \tnorm{\boldsymbol{p}^{s,i} - \boldsymbol{p}_h^{s,i}}_{q,s}^2
      &\le 2\bigg[ (\bar{C}_1\Delta t + \bar{C}_2h^{\ell})^2 + C (\mu^s)^{-1} h^{2\ell} \Delta t \sum_{i=1}^n \norm[0]{ p^{s,i} }_{\ell, \Omega^s}^2 \bigg]
      \\
      &\le 2\bigg[ (\bar{C}_1\Delta t + \bar{C}_2h^{\ell})^2 + C (\mu^s)^{-1} T h^{2\ell} \norm[0]{ p^{s} }_{C^0(J;H^{\ell}(\Omega^s))}^2 \bigg]
      \\
      &\le 2\bigg[ 2\bar{C}_1^2\Delta t^2 + (2\bar{C}_2^2 + C(\mu^s)^{-1}T \norm[0]{ p^{s} }_{C^0(J;H^{\ell}(\Omega^s))}^2)h^{2\ell}\bigg],
    \end{split}
  \end{equation*}
  so that the result follows. 
\end{proof}

\section{Numerical examples}
\label{sec:numerical_examples}

We present some numerical examples using the fully discrete HDG method
\cref{eq:fullydiscreteHDG} to find approximate solutions to the
coupled Stokes and Biot problem
\cref{eq:stokesbiot,eq:interface,eq:bics}. All examples have been
implemented using Netgen/NGSolve \cite{Schoberl:1997,Schoberl:2014}.

\subsection{Stationary test case}
\label{ss:stat_tc_manufactured}

In this first test case we consider the following stationary problem:
\begin{subequations}
  \label{eq:ststokesbiot}
  \begin{align}
    \label{eq:ststokesbiot_1}
    -\nabla \cdot \sigma^j &= f^j && \text{in } \Omega^j, \quad j=s,b,
    \\
    \label{eq:ststokesbiot_2a}
    -\nabla \cdot u^s 
    &= 0 && \text{in } \Omega^s,
    \\
    \label{eq:ststokesbiot_2b}
    -\nabla \cdot u^b + \lambda^{-1}(\alpha p^p - p^b)
    &= 0 && \text{in } \Omega^b,
    \\
    \label{eq:ststokesbiot_3}
    c_0 \tau p^p + \alpha \tau \lambda^{-1} (\alpha p^p - p^b) + \nabla \cdot z
                           &= g^b && \text{in }\Omega^b,
    \\
    \label{eq:ststokesbiot_4}
    \kappa^{-1} z + \nabla p^p
    &= 0 && \text{in }\Omega^b,
  \end{align}  
\end{subequations}
with boundary conditions
\begin{subequations}
  \label{eq:stbics}
  \begin{align}
    \label{eq:stbics-us}
    u^j &= U^j && \text{on } \Gamma^j_D, \quad j=s,b,
    \\
    \label{eq:stbics-sigmas}
    \sigma^j n &= S^j &&  \text{on } \Gamma^j_N, \quad j=s,b,
    \\    
    \label{eq:stbics-pb}
    p^p &= P^p && \text{on } \Gamma^b_P,
    \\
    \label{eq:stbics-zb}
    z \cdot n &= Z && \text{on } \Gamma^b_F,
  \end{align}
\end{subequations}
and interface conditions
\begin{subequations}
  \label{eq:stinterface}
  \begin{align}
    \label{eq:stbc_I_u}
    u^s\cdot n &= \del[1]{\tau u^b + z}\cdot n + M^u & & \text{on } \Gamma_I,
    \\
    \label{eq:stbc_I_ss_sb}
    \sigma^sn &= \sigma^bn + M^s & & \text{on } \Gamma_I,
    \\
    \label{eq:stbc_I_p}
    -(\sigma^s n)\cdot n &= p^p + M^p & & \text{on } \Gamma_I,
    \\    
    \label{eq:stbc_I_slip}
    -2\mu^s\del[0]{\varepsilon(u^s)n}^t &= \gamma (\mu^s/\kappa)^{1/2}(u^s-\partial_t u^b)^t + M^e
                                                     & & \text{on } \Gamma_I.
  \end{align}
\end{subequations}
We consider the unit square domain $\Omega=(0,1)^2$ partitioned as:
$\overline{\Omega} = \overline{\Omega}^s \cup \overline{\Omega}^b$
with $\overline{\Omega}^s = [0,1] \times [1/2,1]$ and
$\overline{\Omega}^b = [0,1] \times [0, 1/2]$. We set
$\Gamma^s_D = \cbr[0]{x\in \Gamma^s:\ x_1=0 \text{ or } x_2=1}$,
$\Gamma^s_N = \cbr[0]{x \in \Gamma^s:\ x_1=1}$,
$\Gamma^b_P = \Gamma^b_D = \cbr[0]{x\in \Gamma^b:\ x_1=0 \text{ or }
  x_2=0}$, and
$\Gamma^b_N = \Gamma^b_F = \cbr[0]{x\in \Gamma^b:\ x_1=1}$. The source
terms $f^s$, $f^b$, and $g^b$, the boundary data $U^s$, $U^b$, $S^s$,
$S^b$, $P^p$, and $Z$, and the interface data $M^u$, $M^s$, $M^p$, and
$M^e$ are chosen such that the exact solution is given by:
\begin{subequations}
  \begin{align*}
    u^s &= 
          \begin{bmatrix}
            \pi x_1\cos(\pi x_1 x_2) + 1 \\
            -\pi x_2 \cos(\pi x_1 x_2) + 2x_1
          \end{bmatrix},
        &&
           p^s = \sin(3x_1)\cos(4x_2), 
    \\
    u^b &=
          \begin{bmatrix}
            \cos(4x_1)\cos(3x_2) \\
            \sin(5x_1)\cos(2x_2)
          \end{bmatrix},
        &&
           p^p = \sin(3x_1x_2),
  \end{align*}
\end{subequations}
Note that $p^b = -\lambda \nabla \cdot u^b + \alpha p^p$ and
$z = -\kappa \nabla p^p$.

We choose the following parameters: $\mu^s = 10^{-2}$,
$\mu^b = 10^{-3}$, $\alpha=0.2$, $\lambda = 10^2$, $\kappa = 10^{-2}$,
$c_0 = 10^{-2}$, $\gamma = 0.3$, $\tau = 10^{-2}$. We choose the
interior penalty parameters as $\beta^s = \beta^b = 8k^2$, where $k$
is the polynomial degree.

We present the errors in the $L^2$-norm and rates of convergence for
all unknowns in
\cref{tab:stat_tc_manufactured_omegaS,tab:sstat_tc_manufactured_omegaB}
for polynomial degrees $k=1$, $k=2$, and $k=3$. We observe optimal
rates of convergence for all unknowns.

\begin{table}
  \small
  \centering {
    \begin{tabular}{cccccc}
      \hline
      Cells
      & $\norm[0]{u_h^s - u^s}_{\Omega^s}$ & Rate & $\norm[0]{p_h^s - p^s}_{\Omega^s}$
      & Rate & $\norm{\nabla \cdot u_h^s}_{\Omega^s}$ \\
      \hline
      \multicolumn{6}{l}{$k=1$} \\
    152 & 1.6e-02 &   - & 4.3e-02 &   - & 2.2e-15 \\ 
    608 & 4.4e-03 & 1.9 & 2.3e-02 & 0.9 & 6.8e-15 \\ 
   2432 & 1.1e-03 & 2.0 & 1.1e-02 & 1.0 & 1.3e-14 \\ 
   9728 & 2.8e-04 & 2.0 & 5.6e-03 & 1.0 & 2.5e-14 \\ 
  38912 & 7.0e-05 & 2.0 & 2.8e-03 & 1.0 & 5.1e-14 \\      
      \multicolumn{6}{l}{$k=2$} \\
    152 & 1.4e-03 &   - & 3.8e-03 &   - & 1.4e-14 \\ 
    608 & 2.0e-04 & 2.8 & 1.1e-03 & 1.8 & 3.5e-14 \\ 
   2432 & 2.5e-05 & 3.0 & 2.7e-04 & 2.0 & 7.1e-14 \\ 
   9728 & 3.2e-06 & 3.0 & 6.6e-05 & 2.0 & 1.4e-13 \\ 
  38912 & 4.0e-07 & 3.0 & 1.7e-05 & 2.0 & 2.8e-13 \\       
      \multicolumn{6}{l}{$k=3$} \\
    152 & 6.2e-05 &   - & 1.9e-04 &   - & 1.8e-13 \\ 
    608 & 4.8e-06 & 3.7 & 3.2e-05 & 2.6 & 3.9e-13 \\ 
   2432 & 3.0e-07 & 4.0 & 3.9e-06 & 3.0 & 7.8e-13 \\ 
   9728 & 1.8e-08 & 4.0 & 4.9e-07 & 3.0 & 1.6e-12 \\ 
  38912 & 1.1e-09 & 4.0 & 6.1e-08 & 3.0 & 3.2e-12 \\       
      \hline
    \end{tabular}
  } \caption{Errors and rates of convergence in $\Omega^s$ for
different polynomial degrees $k$ for the test case described in
\cref{ss:stat_tc_manufactured}.}
  \label{tab:stat_tc_manufactured_omegaS}
\end{table}

\begin{table}
  \tiny
    \centering {
      \begin{tabular}{ccccccccc}
        \hline
        Cells &
        $\norm[0]{u_h^b - u^b}_{\Omega^b}$ & Rate &
        $\norm[0]{p^b_h - p^b}_{\Omega^b}$
        & Rate & $\norm[0]{z_h-z}_{\Omega^b}$ & Rate & $\norm[0]{p_h^p - p^p}_{\Omega^b}$ & Rate \\
        \hline
        \multicolumn{9}{l}{$k=1$} \\
        152 & 2.0e-01 &   - & 2.1e+01 &   - & 7.2e-04 &   - & 2.8e-02 &   -  \\
        608 & 4.5e-02 & 2.1 & 1.2e+01 & 0.9 & 2.5e-04 & 1.5 & 1.5e-02 & 0.9 \\
        2432 & 1.2e-02 & 2.0 & 5.9e+00 & 1.0 & 7.6e-05 & 1.7 & 7.7e-03 & 1.0 \\
        9728 & 2.9e-03 & 2.0 & 2.9e+00 & 1.0 & 1.7e-05 & 2.1 & 3.9e-03 & 1.0 \\
        38912 & 7.2e-04 & 2.0 & 1.5e+00 & 1.0 & 3.4e-06 & 2.4 & 1.9e-03 & 1.0 \\
        \multicolumn{9}{l}{$k=2$} \\
        152 & 5.5e-03 &   - & 1.4e+00 &   - & 4.6e-05 &   - & 1.2e-03 &   -  \\
        608 & 2.5e-04 & 4.4 & 5.8e-01 & 1.3 & 5.6e-06 & 3.0 & 3.9e-04 & 1.6 \\
        2432 & 2.8e-05 & 3.2 & 1.4e-01 & 2.0 & 7.4e-07 & 2.9 & 9.7e-05 & 2.0 \\
        9728 & 3.2e-06 & 3.1 & 3.6e-02 & 2.0 & 7.7e-08 & 3.3 & 2.4e-05 & 2.0 \\
        38912 & 3.9e-07 & 3.1 & 9.0e-03 & 2.0 & 7.2e-09 & 3.4 & 6.0e-06 & 2.0 \\
        \multicolumn{9}{l}{$k=3$} \\
        152 & 1.6e-04 &   - & 7.1e-02 &   - & 3.9e-06 &   - & 3.6e-05 &   -  \\
        608 & 6.4e-06 & 4.7 & 1.6e-02 & 2.2 & 2.5e-07 & 4.0 & 8.1e-06 & 2.2 \\
        2432 & 3.8e-07 & 4.1 & 2.0e-03 & 3.0 & 1.5e-08 & 4.0 & 1.0e-06 & 3.0 \\
        9728 & 2.3e-08 & 4.0 & 2.5e-04 & 3.0 & 7.8e-10 & 4.3 & 1.3e-07 & 3.0 \\
        38912 & 1.5e-09 & 4.0 & 3.1e-05 & 3.0 & 3.6e-11 & 4.4 & 1.6e-08 & 3.0 \\
        \hline
      \end{tabular}      
    }
    \caption{Errors and rates of convergence in $\Omega^b$ for
      different polynomial degrees $k$ for the test case described in
      \cref{ss:stat_tc_manufactured}.}
  \label{tab:sstat_tc_manufactured_omegaB}
\end{table}

\subsection{Time-dependent test case}
\label{ss:td_tc_manufactured}

We now consider the time-dependent problem \cref{eq:stokesbiot} with
boundary and interface conditions given by, respectively,
\begin{subequations}
  \begin{align*}
    u^j &= U^j && \text{on } \Gamma^j_D\times J, \quad j=s,b,
    \\
    \sigma^j n &= S^j &&  \text{on } \Gamma^j_N \times J, \quad j=s,b,
    \\    
    p^p &= P^p && \text{on } \Gamma^b_P\times J,
    \\
    z\cdot n &= Z && \text{on } \Gamma^b_F\times J,
  \end{align*}
\end{subequations}
and
\begin{subequations}
  \begin{align*}
    u^s\cdot n &= \del[1]{\partial_t u^b + z}\cdot n + M^u & & \text{on } \Gamma_I \times J,
    \\
    \sigma^sn &= \sigma^bn + M^s & & \text{on } \Gamma_I \times J,
    \\
    -(\sigma^s n)\cdot n &= p^p + M^p& & \text{on } \Gamma_I \times J,
    \\    
    -2\mu^s\del[0]{\varepsilon(u^s)n}^t &= \gamma (\mu^s/\kappa)^{1/2}(u^s-\partial_t u^b)^t + M^e
                                                     & & \text{on } \Gamma_I \times J.
  \end{align*}
\end{subequations}
We consider the same domain and partitioning of the boundary as in
\cref{ss:stat_tc_manufactured}. The source terms $f^s$, $f^b$, and
$g^b$, the boundary data $U^s$, $U^b$, $S^s$, $S^b$, $P^p$, and $Z$,
and the interface data $M^u$, $M^s$, $M^p$, and $M^e$ are chosen such
that the exact solution is given by:
\begin{subequations}
  \begin{align*}
    u^s &= 
          \begin{bmatrix}
            \pi x_1\cos(\pi (x_1 x_2-t)) + 1 \\
            -\pi x_2 \cos(\pi (x_1 x_2-t)) + 2x_1
          \end{bmatrix},
        &&
           p^s = \sin(3x_1)\cos(4(x_2-t)), 
    \\
    u^b &=
          \begin{bmatrix}
            \sin(10\pi t)\cos(4(x_1-t))\cos(3x_2) \\
            \sin(10\pi t)\sin(5x_1)\cos(2(x_2-t))
          \end{bmatrix},
        &&
           p^p = \sin(3(x_1x_2-t)),
  \end{align*}
\end{subequations}
Note that $p^b = -\lambda \nabla \cdot u^b + \alpha p^p$ and
$z = -\kappa \nabla p^p$.

We choose the following parameters: $\mu^s = 10^{-2}$,
$\mu^b = 10^{-3}$, $\alpha=0.2$, $\lambda = 10^2$, $\kappa = 10^{-2}$,
$c_0 = 10^{-2}$, and $\gamma = 0.3$. We choose the interior penalty
parameters as $\beta^s = \beta^b = 8k^2$, where $k$ is the polynomial
degree. To avoid needing to take very small time steps, we implement
the two-step Backward Differentiation Formulae (BDF2) time-stepping
method. We choose the time step $\Delta t = \tfrac{1}{10} h^{3/2}$ and
consider the time interval $J=[0, 0.01]$.

We present the errors in the $L^2$-norm and rates of convergence for
all unknowns in
\cref{tab:td_tc_manufactured_omegaS,tab:td_tc_manufactured_omegaB} for
polynomial degree $k=2$. We observe optimal rates of convergence for
all unknowns.

\begin{table}
  \small
  \centering {
    \begin{tabular}{cccccc}
      \hline
      Cells
      & $\norm[0]{u_h^s - u^s}_{\Omega^s}$ & Rate & $\norm[0]{p_h^s - p^s}_{\Omega^s}$
      & Rate & $\norm{\nabla \cdot u_h^s}_{\Omega^s}$ \\
      \hline
    152 & 1.5e-03 &   - & 3.8e-03 &   - & 1.4e-14 \\ 
    576 & 1.1e-04 & 3.8 & 8.9e-04 & 2.1 & 2.8e-14 \\ 
   2348 & 1.2e-05 & 3.2 & 2.1e-04 & 2.1 & 5.7e-14 \\ 
   9526 & 1.4e-06 & 3.1 & 5.2e-05 & 2.1 & 1.1e-13 \\       
      \hline
    \end{tabular}
  } \caption{Errors and rates of convergence in $\Omega^s$ for
    polynomial degree $k=2$ for the test case described in
    \cref{ss:td_tc_manufactured}.}
  \label{tab:td_tc_manufactured_omegaS}
\end{table}

\begin{table}
  \tiny
  \centering {
    \begin{tabular}{ccccccccc}
      \hline
      Cells
      & $\norm[0]{u_h^b - u^b}_{\Omega^b}$ & Rate & $\norm[0]{p^b_h - p^b}_{\Omega^b}$
      & Rate & $\norm[0]{z_h-z}_{\Omega^b}$ & Rate & $\norm[0]{p_h^p - p^p}_{\Omega^b}$ & Rate \\
      \hline
    152 & 1.9e-04 &   - & 4.4e-01 &   - & 4.3e-03 &   - & 5.4e-02 &   -  \\ 
    576 & 2.5e-05 & 3.0 & 1.0e-01 & 2.1 & 7.7e-04 & 2.5 & 9.7e-03 & 2.5 \\ 
   2348 & 2.8e-06 & 3.1 & 2.6e-02 & 2.0 & 1.0e-04 & 2.9 & 1.3e-03 & 2.9 \\ 
   9526 & 3.4e-07 & 3.0 & 6.2e-03 & 2.0 & 1.3e-05 & 2.9 & 1.7e-04 & 2.9 \\      
    \hline
    \end{tabular}
  }
  \caption{Errors and rates of convergence in $\Omega^b$ for polynomial degree $k=2$ for the test case described in
    \cref{ss:td_tc_manufactured}.}
  \label{tab:td_tc_manufactured_omegaB}
\end{table}

\subsection{Coupling of surface/subsurface flow}
\label{ss:csursubflow}

In this final example we consider an example proposed in \cite[Section
8.2]{Li:2022}. For this we consider the domain
$\Omega = (0,2) \times (-1,1)$ with $\Omega^s = (0,2) \times (0, 1)$
and $\Omega^b = (0,2) \times (-1, 0)$ and the time interval $J=(0, T)$
with $T=3$. The body forces, source/sink terms, and initial conditions
are set as $f^s=0$, $f^b=0$, $g^b=0$, $p_0=0$, and $u_0=0$. We
consider three parameter sets: (1)
$(\kappa, c_0, \lambda, \mu^b) = (1, 1, 1, 1)$; (2)
$(\kappa, c_0, \lambda, \mu^b) = (10^{-4}, 10^{-4}, 10^6, 1)$; and (3)
$(\kappa, c_0, \lambda, \mu^b) = (10^{-4}, 10^{-4}, 10^6, 10^6)$. The
remaining parameters are chosen as $\mu^s=1$, $\alpha=1$, and
$\gamma=1$.

Let $\Gamma_D^s = \partial\Omega^s \cap \partial\Omega$,
$\Gamma_N^b = \Gamma_P^b = \cbr[0]{x \in \partial\Omega^b\, : \,
  x_2=-1}$, and
$\Gamma_D^b=\Gamma_F^b = \partial\Omega^b \setminus (\Gamma_I \cup
\Gamma_N^b)$. We impose the following boundary conditions:
\begin{align*}
  u^s &= \sbr[1]{-20 x_2(x_2-1)(2-x_1), 0}^T, && \text{on } \Gamma_D^s,
  \\
  u^b &= 0, \quad z\cdot n = 0, && \text{on } \Gamma_D^b = \Gamma_F^b,
  \\
  p^p &= 0, \quad \sigma^bn=0, && \text{on } \Gamma_N^b=\Gamma_P^b.
\end{align*}

We compute the solution on an unstructured simplicial mesh consisting
of 9508 elements, using $k=2$, and a time step of $\Delta t = 0.06$.

\begin{figure}
 \centering
 \subfloat[Velocity. \label{fig:c1_a}]
 {\includegraphics[width=0.31\textwidth]{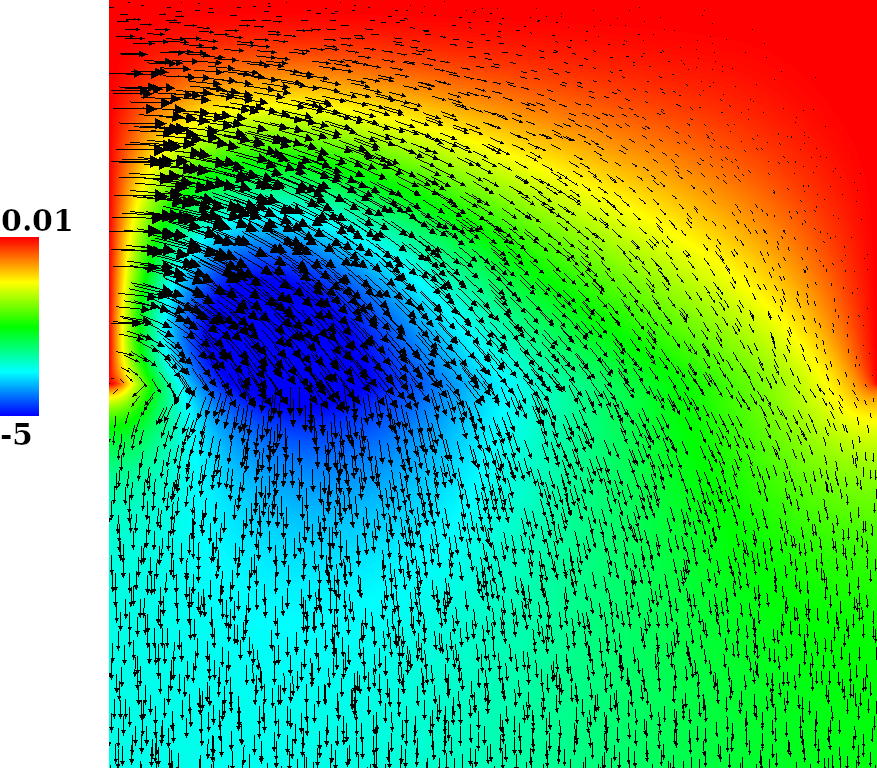}} \quad
 \subfloat[Stresses. \label{fig:c1_b}]
 {\includegraphics[width=0.31\textwidth]{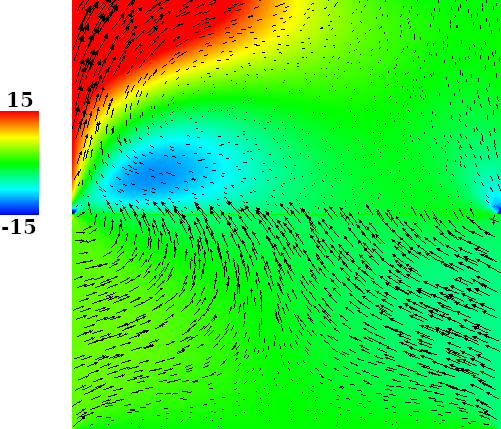}} \quad
 \subfloat[Stresses. \label{fig:c1_c}]
 {\includegraphics[width=0.31\textwidth]{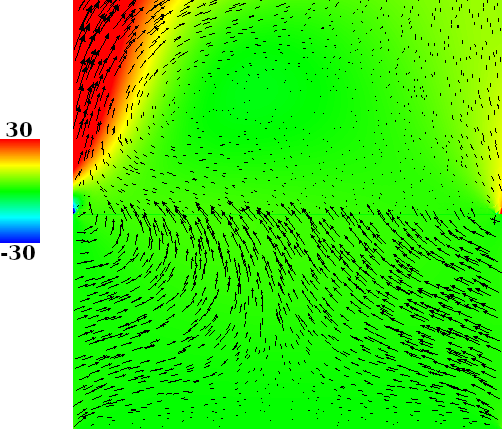}}
 \\
 \subfloat[Displacement. \label{fig:c1_d}]
 {\includegraphics[width=0.48\textwidth]{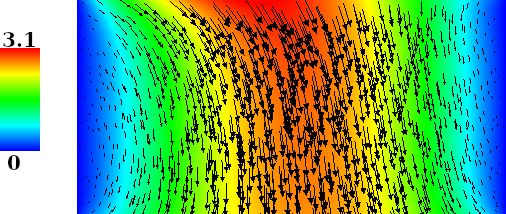}} \quad
 \subfloat[Pore pressure. \label{fig:c1_e}]
 {\includegraphics[width=0.48\textwidth]{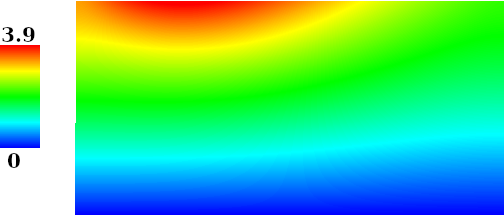}}  
 \caption{Test case from~\cref{ss:csursubflow} with parameter set 1,
    $(\kappa, c_0, \lambda, \mu^b) = (1, 1, 1, 1)$.  Top left: $u^s$
    and $z+\partial_tu^b$ (arrows) and $u^s_2$ and
    $z_2+\partial_tu^b_2$ (color).  Top middle:
    $-(\sigma^s_{12},\sigma_{22}^s)^T$ and
    $-(\sigma^b_{12},\sigma_{22}^b)^T$ (arrows) and $-\sigma_{12}^s$
    and $-\sigma_{12}^b$ (color).  Top right:
    $-(\sigma^s_{12},\sigma_{22}^s)^T$ and
    $-(\sigma^b_{12},\sigma_{22}^b)^T$ (arrows) and $-\sigma_{22}^s$
    and $-\sigma_{22}^b$ (color).  Bottom left: $u^b$ (arrows) and
    $|u^b|$ (color). Bottom right: $p^p$.}
 \label{fig:case1}
\end{figure}

\begin{figure}
 \centering
 \subfloat[Velocity. \label{fig:c2_a}]
 {\includegraphics[width=0.31\textwidth]{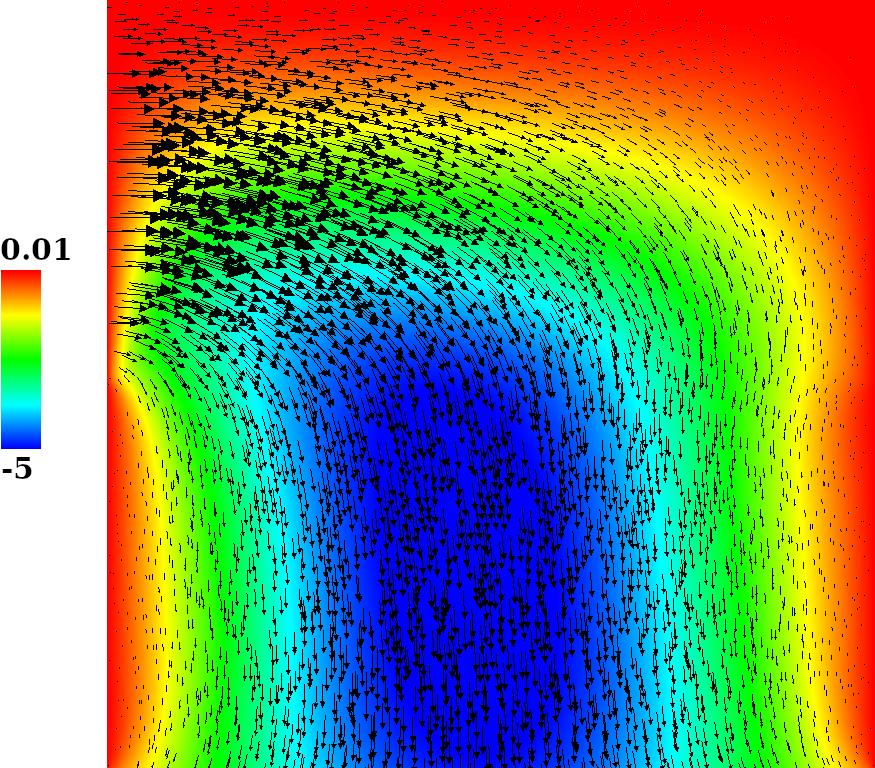}} \quad
 \subfloat[Stresses. \label{fig:c2_b}]
 {\includegraphics[width=0.31\textwidth]{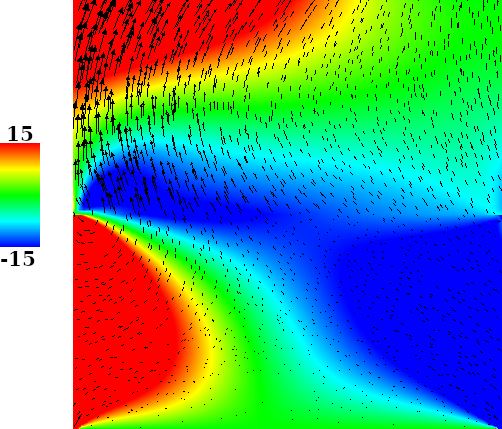}} \quad
 \subfloat[Stresses. \label{fig:c2_c}]
 {\includegraphics[width=0.31\textwidth]{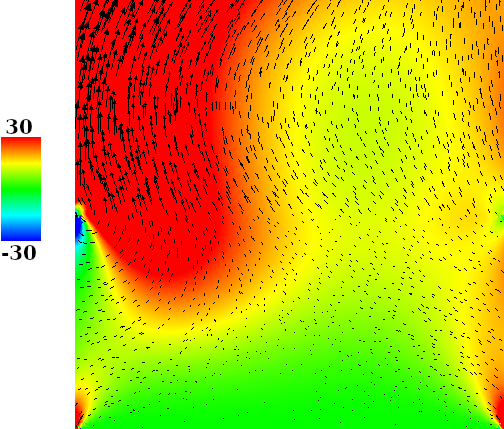}}
 \\
 \subfloat[Displacement. \label{fig:c2_d}]
 {\includegraphics[width=0.48\textwidth]{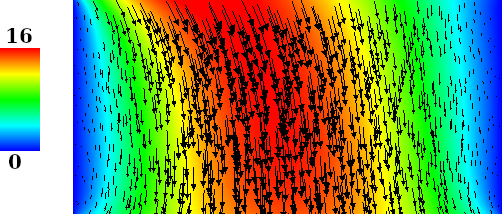}} \quad
 \subfloat[Pore pressure. \label{fig:c2_e}]
 {\includegraphics[width=0.48\textwidth]{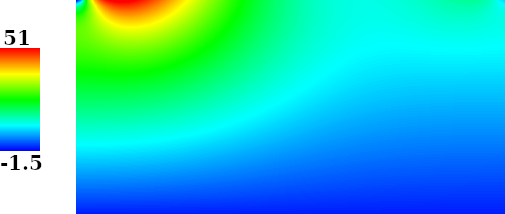}}  
 \caption{Test case from~\cref{ss:csursubflow} with parameter set 2,
    $(\kappa, c_0, \lambda, \mu^b) = (10^{-4}, 10^{-4}, 10^6, 1)$.
    Top left: $u^s$ and $z+\partial_tu^b$ (arrows) and $u^s_2$ and
    $z_2+\partial_tu^b_2$ (color).  Top middle:
    $-(\sigma^s_{12},\sigma_{22}^s)^T$ and
    $-(\sigma^b_{12},\sigma_{22}^b)^T$ (arrows) and $-\sigma_{12}^s$
    and $-\sigma_{12}^b$ (color).  Top right:
    $-(\sigma^s_{12},\sigma_{22}^s)^T$ and
    $-(\sigma^b_{12},\sigma_{22}^b)^T$ (arrows) and $-\sigma_{22}^s$
    and $-\sigma_{22}^b$ (color).  Bottom left: $u^b$ (arrows) and
    $|u^b|$ (color). Bottom right: $p^p$.  }
 \label{fig:case2}
\end{figure}

\begin{figure}
 \centering
 \subfloat[Velocity. \label{fig:c3_a}]
 {\includegraphics[width=0.31\textwidth]{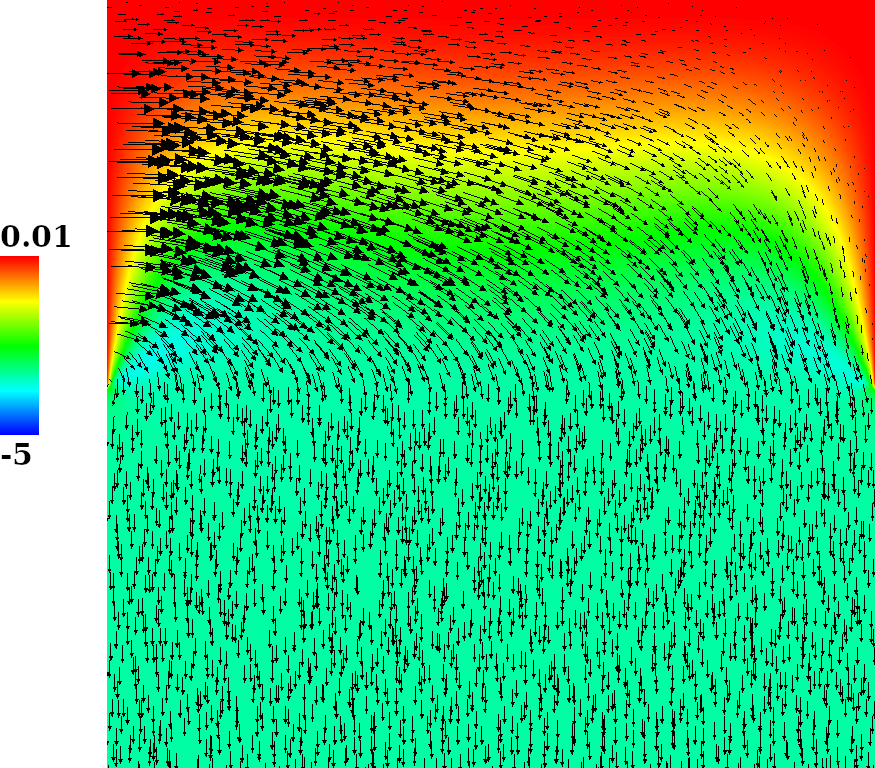}} \quad
 \subfloat[Stresses. \label{fig:c3_b}]
 {\includegraphics[width=0.31\textwidth]{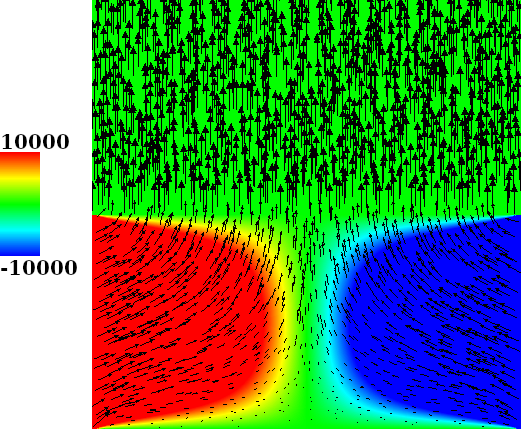}} \quad
 \subfloat[Stresses. \label{fig:c3_c}]
 {\includegraphics[width=0.31\textwidth]{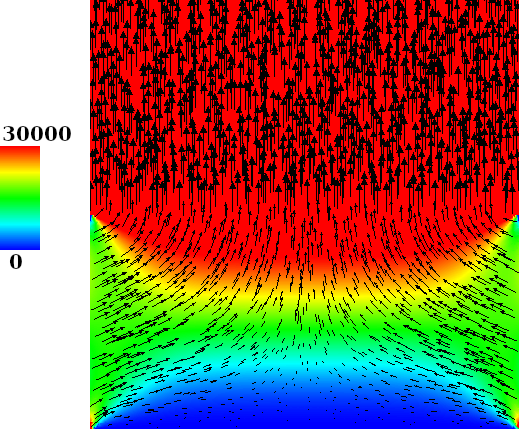}}
 \\
 \subfloat[Displacement. \label{fig:c3_d}]
 {\includegraphics[width=0.48\textwidth]{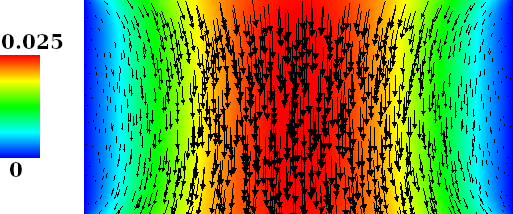}} \quad
 \subfloat[Pore pressure. \label{fig:c3_e}]
 {\includegraphics[width=0.48\textwidth]{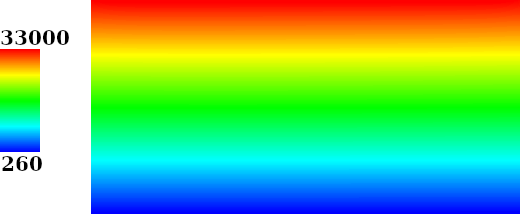}}  
 \caption{Test case from~\cref{ss:csursubflow} with parameter set 3,
    $(\kappa, c_0, \lambda, \mu^b) = (10^{-4}, 10^{-4}, 10^6, 10^6)$.
    Top left: $u^s$ and $z+\partial_tu^b$ (arrows) and $u^s_2$ and
    $z_2+\partial_tu^b_2$ (color).  Top middle:
    $-(\sigma^s_{12},\sigma_{22}^s)^T$ and
    $-(\sigma^b_{12},\sigma_{22}^b)^T$ (arrows) and $-\sigma_{12}^s$
    and $-\sigma_{12}^b$ (color).  Top right:
    $-(\sigma^s_{12},\sigma_{22}^s)^T$ and
    $-(\sigma^b_{12},\sigma_{22}^b)^T$ (arrows) and $-\sigma_{22}^s$
    and $-\sigma_{22}^b$ (color).  Bottom left: $u^b$ (arrows) and
    $|u^b|$ (color). Bottom right: $p^p$.  }
 \label{fig:case3}
\end{figure}

We plot the solution obtained with the three different parameter sets
in \cref{fig:case1,fig:case2,fig:case3}. The results compare well to
those obtained by the locking-free method of \cite{Li:2022}; the
solution does not exhibit locking or oscillations despite Poisson
ratio $\nu=0.4999995$ (for parameter set 2) and despite modeling a
very stiff poroelastic medium (parameter set 3). Furthermore, we
observe from \cref{fig:c1_a,fig:c2_a,fig:c3_a} that the second
component of the velocity is continuous across the interface, i.e.,
mass is conserved at the interface. From
\cref{fig:c1_b,fig:c2_b,fig:c3_b} and from
\cref{fig:c1_c,fig:c2_c,fig:c3_c} we observe that
$-\sigma^s_{12}=-\sigma^b_{12}$ and $-\sigma_{22}^s=-\sigma_{22}^b$
implying conservation of momentum on the interface.

\section{Conclusions}
\label{sec:conclusions}

We introduced an HDG method for the coupled Stokes--Biot problem that
is provably robust in the incompressible limit, $\lambda \to \infty$
and $c_0 \to 0$. Consistency was shown for the semi-discrete case
while well-posedness and a priori error estimates were determined
after combining the HDG method with backward Euler
time-stepping. Furthermore, we showed that the discrete velocities and
displacement are $H(\text{div})$-conforming and that the
compressibility equations are satisfied pointwise by the numerical
solution on the elements. Mass is conserved pointwise on the elements
for the semi-discrete problem (up to the error of the $L^2$-projection
of the source/sink term into the discrete pore pressure space).
Finally, numerical examples demonstrate optimal rates of convergence
for all unknowns in the $L^2$-norm and that the numerical method is
locking-free.

\subsubsection*{Acknowledgements}

Aycil Cesmelioglu and Jeonghun J. Lee gratefully acknowledge support by the National Science
Foundation (grant numbers DMS-2110782 and DMS-2110781) and Sander Rhebergen
gratefully acknowledges support from the Natural Sciences and
Engineering Research Council of Canada through the Discovery Grant
program (RGPIN-05606-2015).

\bibliography{references}
\bibliographystyle{plain}
\appendix
\section{Proof of the inf-sup condition \cref{eq:total-pressure-inf-sup-alt}}
\label{ap:infsupproof}

The inf-sup condition was proven for the case of homogeneous Dirichlet
boundary conditions on the whole boundary of the domain in \cite[Lemma
4.4]{Rhebergen:2017} and \cite[Lemma 1]{Rhebergen:2018b} for the HDG
method and \cite[Lemma 8]{Rhebergen:2020} for a variation of the HDG
method. Here we generalize these proofs to the case where homogeneous
Dirichlet boundary conditions are only posed on part of the domain
boundary. The proof proceeds in three steps.

\textbf{Step 1.} Let $\Pi_V^j:[H^1(\Omega^j)]^d \to V_h^j$ be the BDM
interpolation operator \cite[Section III.3]{Brezzi:book} and define
the norm $\norm[0]{v_h^j}_{1,h,j}$ as
\begin{equation*}
  \norm[0]{v_h^j}_{1,h,j}^2 := \norm[0]{\varepsilon(v_h^j)}_{\Omega^j}^2 + \sum_{F \in \mathcal{F}^j} h_F^{-1} \norm[0]{\jump{v_h^j}}_F^2,
\end{equation*}
where $h_F$ is the radius of the face $F$. Note that
\begin{equation*}
  \tnorm{(v_h, \av{v_h})}_{v,j}^2
  = \sum_{K \in \mathcal{T}^j}\del[1]{ \norm[0]{\varepsilon(v_h)}_K^2 + h_K^{-1}\norm[0]{v_h - \av{v_h}}_{\partial K}^2}.
\end{equation*}
On boundary facets $\av{v_h}=v_h$ and so
\begin{equation*}
  \begin{split}
    \sum_{K \in \mathcal{T}^j} h_K^{-1}\norm[0]{(v_h - \av{v_h})}_{\partial K}^2
    =& \sum_{F \in \mathcal{F}^j_{int}} \del[1]{h_{K^+}^{-1}\norm[0]{(v_h^+ - \av{v_h})}_{F}^2 + h_{K^-}^{-1}\norm[0]{(v_h^- - \av{v_h})}_{F}^2}
    \\
    =& \sum_{F \in \mathcal{F}^j_{int}} \del[1]{\tfrac{1}{4}h_{K^+}^{-1}\norm[0]{v_h^+ - v_h^-}_{F}^2 + \tfrac{1}{4}h_{K^-}^{-1}\norm[0]{v_h^- - v_h^+}_{F}^2}
    \\
    =& \sum_{F \in \mathcal{F}^j_{int}} \tfrac{1}{4}\del[0]{h_{K^+}^{-1} + h_{K^-}^{-1}}\norm[0]{v_h^+ - v_h^-}_{F}^2
    \\
    \le& C\sum_{F \in \mathcal{F}^j_{int}} h_F^{-1}\norm[0]{\jump{v_h}}_{F}^2,
  \end{split}
\end{equation*}
where, assuming shape regularity of the mesh, the inequality is by
equivalence of $h_F$, $h_{K^+}$, and $h_{K^-}$ where $F$ is a facet
shared by $K^+$ and $K^-$. Throughout this proof $C>0$ is a generic
constant independent of $h$. This shows that
$\tnorm{(v_h, \av{v_h})}_{v,j} \le C \norm[0]{v_h^j}_{1,h,j}$. Then,
for all $v \in [H^1(\Omega^j)]^d$,
\begin{equation}
  \label{eq:aPiVjbound1}
  \tnorm{(\Pi_V^jv, \av{\Pi_V^jv})}_{v,j}
  \le C
  \norm[0]{\Pi_V^jv}_{1,h,j}
  \le C \norm[0]{v}_{1,\Omega^j}^2,
\end{equation}
where the second inequality was shown in the proof of
\cite[Proposition 10]{Hansbo:2002}. We next define
\begin{equation*}
  \begin{split}
    [H^1_{ID,j}(\Omega^j)]^d
    &:= \cbr[1]{v \in [H^1(\Omega^j)]^d\ : \ v=0 \text{ on } \Gamma_I \cup \Gamma_D^j},
    \\
    \tilde{\tilde{\boldsymbol{V}}}_h^j
    &:= \cbr[1]{\boldsymbol{v}_h \in
      \widetilde{\boldsymbol{V}}_h^j\ :\ v_h\in H(\text{div};\Omega^j),\ \bar{v}_h^j
      \cdot n = v_h^j \cdot n \text{ on } \partial \Omega^j}.    
  \end{split}
\end{equation*}
Given $|\Gamma_N^j| > 0$, it was shown in \cite[Lemma
B.1]{Baerland:2017} that there exists a constant $C > 0$ such that for
all $q_h \in Q_h^j$ there is a $v_{q_h} \in [H^1_{ID,j}(\Omega^j)]^d$
that satisfies
\begin{equation}
  \label{eq:avqhboundqh}
  -\nabla \cdot v_{q_h} = q_h \qquad C \norm[0]{v_{q_h}}_{1,\Omega^j} \le \norm{q_h}_{\Omega^j}.
\end{equation}
By \cref{eq:aPiVjbound1,eq:avqhboundqh} we note that
\begin{equation*}
  \tnorm{(\Pi_V^jv_{q_h}, \av{\Pi_V^jv_{q_h}})}_{v,j}
  \le C \norm[0]{v_{q_h}}_{1,\Omega^j}
  \le C \norm{q_h}_{\Omega^j}.
\end{equation*}
Note also that
$(\Pi_V^jv_{q_h}, \av{\Pi_V^jv_{q_h}}) \in
\tilde{\tilde{\boldsymbol{V}}}_h^j$. 
We therefore find that
\begin{equation*}
  \sup_{0 \ne \boldsymbol{v}_h \in \tilde{\tilde{\boldsymbol{V}}}_h^j}\frac{b_h^j((q_h^j,0),\boldsymbol{v}_h)}{\tnorm{\boldsymbol{v}_h}_{v,j}}
  \ge
  \frac{b_h^j((q_h^j,0),(\Pi_V^jv_{q_h}, \av{\Pi_V^jv_{q_h}}))}{\tnorm{(\Pi_V^jv_{q_h}, \av{\Pi_V^jv_{q_h}})}_{v,j}}
  \ge
  C \norm[0]{q_h^j}_{\Omega^j}.
\end{equation*}

\textbf{Step 2.} Noting that
$(v_h, 0) \in \widetilde{\boldsymbol{V}}_h^j$, there
exists a $C$ such that
\begin{equation}
  \sup_{\boldsymbol{0}\ne \boldsymbol{v}_h \in \widetilde{\boldsymbol{V}}_h^j}
  \frac{b_h^j((0,\bar{q}_h),\boldsymbol{v}_h)}{\tnorm{\boldsymbol{v}_h}_{v,j}} \ge
  \sup_{0\ne v_h \in V_h^j} \frac{b_h^j((0,\bar{q}_h),(v_h,0))}{\tnorm{(v_h,0)}_{v,j}} \ge C \tnorm{(0, \bar{q}_h)}_{q,j} \quad
  \forall \bar{q}_h \in \bar{Q}_h^j,
\end{equation}
where the second inequality was shown in the proof of \cite[Lemma
3]{Rhebergen:2018b}. (Although \cite[Lemma 3]{Rhebergen:2018b} assumed
quasi-conformity, the result can be extended to shape-regular meshes.)

\textbf{Step 3.} Define
$b_1^j((q_h^j,0),\boldsymbol{v}_h) =
b_h^j((q_h^j,0),\boldsymbol{v}_h)$ and
$b_2^j((0, \bar{q}_h), \boldsymbol{v}_h)=b_h^j((0, \bar{q}_h),
\boldsymbol{v}_h)$. Noting that
\begin{equation*}
  \tilde{\tilde{\boldsymbol{V}}}_h^j
  =
  \cbr[1]{\boldsymbol{v}_h \in \widetilde{\boldsymbol{V}}_h^j \ : \ b_2^j((0,\bar{q}_h), \boldsymbol{v}_h) = 0 \ \forall \bar{q}_h \in \bar{Q}_h^j},
\end{equation*}
the result follows after applying \cite[Theorem 3.1]{Howell:2011}.

\section{Proof of the inf-sup condition \cref{eq:inf-sup}}
\label{ap:infsupprooffull}

The proof is similar to that given in \cref{ap:infsupproof}. It is
given here for completeness. The proof again proceeds in three steps.

\textbf{Step 1.} Let $\Pi_V:[H^1(\Omega)]^d \to V_h$ be the BDM
interpolation operator \cite[Section III.3]{Brezzi:book} and define
the norm $\norm[0]{v_h}_{1,h}$ as
\begin{equation*}
  \norm[0]{v_h}_{1,h}^2 := \norm[0]{\varepsilon(v_h)}_{\Omega}^2 + \sum_{F \in \mathcal{F}} h_F^{-1} \norm[0]{\jump{v_h}}_F^2,
\end{equation*}
where $h_F$ is the radius of the face $F$. As in \cref{ap:infsupproof}
we have that for all $v \in [H^1(\Omega)]^d$,
\begin{equation}
  \label{eq:aPiVjbound1f}
  \tnorm{(\Pi_Vv, \av{\Pi_Vv}, \av{\Pi_Vv})}_{v}
  \le C
  \norm[0]{\Pi_Vv}_{1,h}
  \le C \norm[0]{v}_{1,\Omega},
\end{equation}
where, in this proof, $C>0$ is a generic constant independent of
$h$. Let $\Gamma_D = \Gamma_D^s \cup \Gamma_D^b$ and
$\Gamma_N = \Gamma_N^s \cup \Gamma_N^b$ and define
\begin{equation*}
  \begin{split}
    [H^1_{D}(\Omega)]^d
    &:= \cbr[1]{v \in [H^1(\Omega)]^d\ : \ v=0 \text{ on } \Gamma_D},
    \\
    \tilde{\tilde{\boldsymbol{V}}}_h
    &:= \cbr[1]{\boldsymbol{v}_h \in
      \widehat{\boldsymbol{V}}_h\ :\ v_h\in H(\text{div};\Omega),\ \bar{v}_h^s \cdot n = \bar{v}_h^b \cdot n = v_h \cdot n \text{ on } \partial \Omega \cup \Gamma_I}.
  \end{split}
\end{equation*}
Given $|\Gamma_N| > 0$ there exists a constant $C > 0$ such that
for all $q_h \in Q_h$ there is a $v_{q_h} \in [H^1_{D}(\Omega)]^d$
that satisfies (see \cite[Lemma B.1]{Baerland:2017})
\begin{equation}
  \label{eq:avqhboundqhf}
  -\nabla \cdot v_{q_h} = q_h \qquad C \norm[0]{v_{q_h}}_{1,\Omega} \le \norm{q_h}_{\Omega}.
\end{equation}
By \cref{eq:aPiVjbound1f,eq:avqhboundqhf},
\begin{equation*}
  \tnorm{(\Pi_Vv_{q_h}, \av{\Pi_Vv_{q_h}}, \av{\Pi_Vv_{q_h}})}_{v}
  \le C \norm[0]{v_{q_h}}_{1,\Omega}
  \le C \norm{q_h}_{\Omega}.
\end{equation*}
Since
$(\Pi_Vv_{q_h}, \av{\Pi_Vv_{q_h}}, \av{\Pi_Vv_{q_h}}) \in
\tilde{\tilde{\boldsymbol{V}}}_h$ we obtain
\begin{equation*}
  \sup_{0 \ne \boldsymbol{v}_h \in \tilde{\tilde{\boldsymbol{V}}}_h}\frac{b_h((q_h,0),\boldsymbol{v}_h)}{\tnorm{\boldsymbol{v}_h}_{v}}
  \ge
  \frac{b_h((q_h,0),(\Pi_Vv_{q_h}, \av{\Pi_Vv_{q_h}}, \av{\Pi_Vv_{q_h}}))}{\tnorm{(\Pi_Vv_{q_h}, \av{\Pi_Vv_{q_h}}, \av{\Pi_Vv_{q_h}})}_{v}}
  \ge
  C \norm[0]{q_h}_{\Omega}.
\end{equation*}

\textbf{Step 2.} Noting that
$(v_h, 0, 0) \in \widehat{\boldsymbol{V}}_h$, there exists a $C$
such that
\begin{equation}
  \sup_{\boldsymbol{0}\ne \boldsymbol{v}_h \in \widehat{\boldsymbol{V}}_h}
  \frac{b_h((0,\bar{q}_h),\boldsymbol{v}_h)}{\tnorm{\boldsymbol{v}_h}_{v}} \ge
  \sup_{0\ne v_h \in V_h} \frac{b_h((0,\bar{q}_h),(v_h,0,0))}{\tnorm{(v_h,0,0)}_{v}} \ge C \tnorm{(0, \bar{q}_h)}_{q} \quad
  \forall \bar{q}_h \in \bar{Q}_h,
\end{equation}
where the second inequality was shown in the proof of \cite[Lemma
3]{Rhebergen:2018b}. 

\textbf{Step 3.} Define
$b_1((q_h,0),\boldsymbol{v}_h) = b_h((q_h,0),\boldsymbol{v}_h)$ and
$b_2((0, \bar{q}_h), \boldsymbol{v}_h)=b_h((0, \bar{q}_h),
\boldsymbol{v}_h)$. Noting that
\begin{equation*}
  \begin{split}
    \tilde{\tilde{\boldsymbol{V}}}_h 
    &=
    \cbr[1]{\boldsymbol{v}_h \in
      \widehat{\boldsymbol{V}}_h\ :\ v_h^j \in H(\text{div};\Omega^j),\ \bar{v}_h^j \cdot n = v_h^j \cdot n \text{ on } \partial \Omega \cup \Gamma_I, \ j=s,b}.
    \\
    &=\cbr[1]{\boldsymbol{v}_h \in \widehat{\boldsymbol{V}}_h \ : \ b_2((0,\bar{q}_h), \boldsymbol{v}_h) = 0 \ \forall \bar{q}_h \in \bar{Q}_h},
  \end{split}
\end{equation*}
the result follows after applying \cite[Theorem 3.1]{Howell:2011}.

\end{document}